\documentclass[a4paper, 10 pt, reqno]{article}
\usepackage{amsfonts, amssymb, amsmath, eucal, amsthm}

\usepackage[all]{xy}

\newtheorem{lemma}{Lemma}
\newtheorem{theorem}[lemma]{Theorem}
\newtheorem*{theorem*}{Theorem}
\newtheorem*{proposition*}{Proposition}
\newtheorem{proposition}[lemma]{Proposition}
\newtheorem{corollary}[lemma]{Corollary}

\theoremstyle{definition}
\newtheorem{definition}[lemma]{Definition}
\newtheorem*{definition*}{Definition}
\newtheorem{note}[lemma]{Note}
\newtheorem{example}[lemma]{Example}
\newtheorem{remark}[lemma]{Remark}
\newtheorem{question}[lemma]{Question}
\newtheorem*{acknowledgments}{Acknowledgments}

\numberwithin{lemma}{section}

\newcommand{\B}{\mathbb{B}}
\newcommand{\R}{\mathbb{R}}

\newcommand{\X}{\mathfrak{X}}
\newcommand{\Y}{\mathfrak{Y}}

\newcommand{\bigddt}{{\frac{\partial}{\partial t}}}
\newcommand{\smallddt}{{\partial/\partial t}}

\newcommand{\Mor}{\mathrm{Mor}}
\newcommand{\Diff}{\mathrm{Diff}}
\newcommand{\StDiff}{\mathrm{StDiff}}
\newcommand{\Gpd}{\mathrm{Gpd}}
\newcommand{\gpd}{\mathrm{Lie}}
\newcommand{\LieGpd}{\mathrm{LieGpd}}
\newcommand{\cone}{\mathrm{Cone}}
\newcommand{\st}{\mathrm{st}}
\newcommand{\Cat}{\mathrm{Cat}}
\newcommand{\Vect}{\mathrm{Vect}}

\DeclareMathOperator{\colim}{\mathrm{colim}}
\newcommand{\Id}{\mathrm{Id}}

\newcommand{\BaezCrans}{MR2068522}
\newcommand{\BaezDolan}{MR1664990} 

\newcommand{\Borceux}{MR1291599}

\newcommand{\BX}{BehrendXu} 
\newcommand{\GH}{GepnerHenriques}  
\newcommand{\JeffJohannes}{JeffJohannes}
\newcommand{\KobayashiNomizuOne}{MR1393940}
\newcommand{\Lerman}{Lerman}
\newcommand{\Heinloth}{Heinloth}
\newcommand{\Milnor}{MR0163331}
\newcommand{\MoerdijkMacLane}{MR1300636}
\newcommand{\Moerdijk}{MR1950948}
\newcommand{\Morse}{MorseInequalities} 

\newcommand{\Zung}{MR2292634} 

\title{Vector Fields and Flows on Differentiable Stacks}
\author{Richard Hepworth\footnote{The author is supported by E.P.S.R.C.~Postdoctoral Research Fellowship EP/D066980.}\\  Department of Pure Mathematics\\ University of Sheffield}
\begin{document}
\maketitle
\begin{abstract}
This paper introduces the notions of vector field and flow on a general differentiable stack.  Our main theorem states that the flow of a vector field on a compact proper differentiable stack exists and is unique up to a uniquely determined 2-cell.  This extends the usual result on the existence and uniqueness of flows on a manifold as well as the author's existing results for orbifolds.  It sets the scene for a discussion of Morse Theory on a general proper stack and also paves the way for the categorification of other key aspects of differential geometry such as the tangent bundle and the Lie algebra of vector fields.
\end{abstract}

\section{Introduction}\label{Introduction}

This paper extends the notions of \emph{vector field} and \emph{flow} from manifolds to differentiable stacks.  It is part of a programme to establish Morse Theory for stacks, where the principal tool will be the negative gradient flow of an appropriate Morse function.  The Morse Inequalities, Morse Homology Theorem and handlebody decompositions are powerful computational and conceptual consequences of Morse Theory that we hope to bring to bear on the study of differentiable stacks, or equivalently, the study Lie groupoids.  The author has already established the Morse Inequalities for orbifolds, which are the proper \'etale differentiable stacks \cite{\Morse}.
 
Our results are an example of \emph{categorification} \cite{\BaezDolan}.  In one sense categorification means taking a familiar structure defined by sets, functions and equations among the functions, and then considering an analogous structure determined by categories, functors and natural isomorphisms among the functors.  More generally, categorification can refer to the process of taking notions phrased inside a $1$-category and establishing analogues inside a higher category; the sense we mentioned first promotes notions from the $1$-category of sets to the $2$-category of categories.  Differentiable stacks, or rather an appropriate subclass like the Deligne-Mumford stacks or proper stacks, are a categorification of manifolds, just as groupoids are a categorification of sets.  What this paper achieves, then, is a categorification of vector fields and flows.  We hope that it will open up the possibility of categorifying other aspects of differential geometry via stacks and, perhaps more interestingly, seeing which categorified structures will appear in the process.  We shall elaborate on this point later.

The paper begins by defining a \emph{tangent stack functor}.  This is a lax functor from the $2$-category of differentiable stacks to itself and extends the functor that sends a manifold to its tangent bundle and a map to its derivative.  This allows us to give our first definition:

\begin{definition*}
A \emph{vector field} on a differentiable stack $\X$ is a pair $(X,a_X)$ consisting of a morphism
\[X\colon\X\to T\X\]
and a $2$-cell
\begin{equation}\label{IntroVectorFieldEquation}\xymatrix{
\X\ar[r]^X\ar@/_6ex/[rr]_{\Id_\X}^{}="1"& T\X\ar[r]^{\pi_\X}\ar@{=>}"1"^{a_X}& \X.
}\end{equation}
Here $\pi_\X\colon T\X\to\X$ is the natural projection map.
\end{definition*}

When $\X$ is a manifold $M$, there are no nontrivial $2$-cells between maps $M\to M$.  Two maps are either equal or are not related by any $2$-cell.  Thus \eqref{IntroVectorFieldEquation} becomes the familiar equation $\pi_M\circ X=\Id_M$ and we recover the usual definition of vector field on $M$.  However for a general stack the equation $\pi_\X\circ X=\Id_\X$ may fail to hold while many different $2$-morphisms $a_X$ exist.  The definition above is typical of categorification: the familiar equation $\pi_M\circ X=\Id_M$ is `weakened' to become the isomorphism $a_X$.  Another prominent feature of categorification is that the isomorphisms by which one weakened the original equations are often subjected to new equations of their own.  This is apparent in the next definition.

\begin{definition*}
Let $X$ be a vector field on $\X$.  A \emph{flow of $X$} is a morphism
\[\Phi\colon\X\times\R\to\X\]
equipped with $2$-cells
\begin{equation}\label{IntroductionTPhiEquation}\xymatrix{
T(\X\times\R)\ar[r]^-{T\Phi}_{}="2"& T\X\\
\X\times\R\ar[u]^{\bigddt}\ar[r]_-\Phi &\X.\ar[u]_X^{}="1"\ar@{=>}^{t_\Phi}"1";"2"
}\end{equation}
and
\begin{equation}\label{IntroductionEPhiEquation}\xymatrix{
\X\ar@/_6ex/[rr]^{}="1"_{\Id}\ar@{^{(}->}[r]^-{t=0} &\X\times\R\ar[r]^-\Phi\ar@{=>}"1"^{e_\Phi} &\X
}\end{equation}
for which the composition of $2$-cells in
\begin{equation}\label{IntroductionCategorifiedDiagram}\xymatrix{
\X\times\R\ar@{<-}@/_12ex/[dd]^{}="1"_{\Id}\ar[r]_{}="3"^-\Phi & \X\ar@{<-}@/^12ex/[dd]_{}="2"^{\Id} \\
T(\X\times\R)\ar[u]\ar[r]^-{T\Phi}_{}="4"\ar@{=>}"1"^-{a_\smallddt} & T\X\ar@{=>}"2"_-{a_X}\ar[u]_{}="5"\\
\X\times\R\ar[u]\ar[r]_\Phi & \X\ar[u]_{}="6"\ar@{=>}"6";"4"^{t_\Phi}\ar@{=>}"5";"3"
}\end{equation}
is trivial.  (The upper square is obtained from the naturality of the projection maps $T\X\to\X$, $T(\X\times\R)\to\X\times\R$.)
\end{definition*}

Consider again the case where $\X$ is a manifold $M$.  Then \eqref{IntroductionTPhiEquation} and \eqref{IntroductionEPhiEquation} become the familiar equations $\partial\Phi/\partial t=X\circ\Phi$ and $\Phi(x,0)=x$ that define the flow of $X$, while the condition on the diagram \eqref{IntroductionCategorifiedDiagram} is vacuous.  In general, though, there may be a choice of $t_\Phi$ and $e_\Phi$, and not all choices of $t_\Phi$ will satisfy the condition \eqref{IntroductionCategorifiedDiagram}.  Again this is typical categorification: familiar equations are weakened to isomorphisms and a new equation is imposed on these isomorphisms.  With this definition we are able to prove the following theorem, which extends the usual result on the existence and uniqueness of flows on manifolds.

\begin{theorem*}
Let $X$ be a vector field on a proper differentiable stack $\X$.
\begin{enumerate}
\item If $X$ has compact support then a flow $\Phi\colon\X\times\R\to\X$ exists.
\item Any two flows
\[\Phi,\Psi\colon\X\times\R\to\X\]
of $X$ are related by a $2$-morphism $\Phi\Rightarrow\Psi$ that is uniquely determined by $e_\Phi$, $e_\Psi$, $t_\Phi$ and $t_\Psi$.
\end{enumerate}
(Recall that $\X$ is \emph{proper} if the diagonal map $\Delta\colon\X\to\X\times\X$ is proper.  This is the case for all manifolds, orbifolds, $S^1$-gerbes, and global quotients by compact Lie groups.)
\end{theorem*}

Where do these results lead?  The definitions and theorems described above ignored some of the finer structures available in the theory of tangent bundles, vector fields and flows:
\begin{itemize}
\item The tangent bundle of a manifold is not just a manifold but a vector bundle.  \item The set of vector fields on a manifold is not just a set but a Lie algebra.
\item The set of vector fields on a compact manifold is isomorphic (by taking flows) to the set of $1$-parameter families of diffeomorphisms.
\end{itemize}
Work in progress builds on the present paper and shows that each of the above statements can be extended to stacks:
\begin{itemize}
\item The tangent stack of a differentiable stack is a \emph{bundle of $2$-vector spaces} on the stack.
\item The groupoid of vector fields on a stack is a \emph{Lie $2$-algebra}.
\item The groupoid of vector fields on a compact proper differentiable stack is equivalent to the groupoid of weak actions of $\R$ on the stack.
\end{itemize}
The $2$-vector spaces and Lie $2$-algebras just mentioned should be understood in the sense of Baez and Crans \cite{\BaezCrans}.  By regarding the tangent stack as a $2$-vector bundle we will be able to consider \emph{Riemannian metrics} on a differentiable stack and so to construct \emph{gradient vector fields}.  The gradient vector field of a Morse function, or rather the flow of the gradient, is the fundamental tool in Morse Theory.

The paper is organized as follows.  In \S\ref{TangentStacksSection} we establish the existence of a tangent stack functor $T\colon\StDiff\to\StDiff$ from stacks on $\Diff$ to stacks on $\Diff$.  This is a lax functor that extends the usual tangent functor given by sending a manifold to its tangent bundle and a map to its derivative.  In \S\ref{VectorFieldsSection} we give the full definition of vector fields and  equivalences of vector fields on a stack.  Several key technical results are proved.  We also define vector fields on a Lie groupoid and prove that these are equivalent to vector fields on the stack of torsors.  In \S\ref{IntegralsFlowsSection} we define integral morphisms and integral $2$-morphisms --- these are the analogues of integral curves in a manifold --- and we give the full definition of flows.  Then we state and prove theorems on the existence, uniqueness and representability of integral morphisms and flows, including the theorem stated in this introduction.  \S\ref{QuotientsSection} explores these results in the case of a global quotient stack $[M/G]$ with $G$ a compact Lie group.  The vector fields on $[M/G]$ are described entirely in terms of $G$-equivariant vector fields on $M$, and their flows are described using the flows of these $G$-equivariant fields.  \S\ref{EtaleSection} explores the results for \'etale stacks, and describes how the present results include as a special case the results proved in \cite{\Morse}.  Finally two appendices recall the fundamental Dictionary Lemma and various properties of proper stacks.

\begin{acknowledgments}
Thanks to David Gepner and Jeff Giansiracusa for many interesting and useful discussions about stacks.  The author is supported by an E.P.S.R.C.~Postdoctoral Research Fellowship, grant number EP/D066980.
\end{acknowledgments}

\setcounter{tocdepth}{1}
\tableofcontents

\section{Tangent Stacks}\label{TangentStacksSection}
Let $\Diff$ denote the category of smooth manifolds and smooth maps, equipped with the usual Grothendieck (pre)topology determined by open coverings.  Then \emph{stacks on $\Diff$}, which are the lax sheaves of groupoids on $\Diff$, form a strict $2$-category that we denote by $\StDiff$.  There is a \emph{Yoneda embedding} $y\colon\Diff\to\StDiff$, and so we can think of stacks on $\Diff$ as a generalization of manifolds.  For readable introductions to the language of differentiable stacks we recommend \cite{\Heinloth}, \cite{\BX}.

Taking tangent bundles and derivatives determines a functor
\[T\colon\Diff\to\Diff\]
that we call the \emph{tangent functor}.  Functoriality of $T$ is nothing but the chain rule.  The projections $\pi_X\colon TX\to X$ together constitute a natural \emph{projection map} $\pi\colon T\Rightarrow\Id$.  The object of this section is to define the `tangent stack' of any stack on $\Diff$ in a functorial way that extends the usual notion of tangent bundle for manifolds.

In \S\ref{TangentConstructionSubsection} we define the lax \emph{tangent stack functor}
\[T^\st\colon\StDiff\to\StDiff,\]
and a lax natural morphism $\pi^\st\colon T^\st\Rightarrow\Id$ called the \emph{projection map}.  In \S\ref{YonedaSubsection} we will show that the functor $T^\st$  satisfies $T^\st\circ y=y\circ T$, so that when restricted to manifolds $T^\st$ is just the usual tangent functor $T\colon\Diff\to\Diff$.

Among all stacks on $\Diff$ it is common to concentrate on the \emph{differentiable stacks}.  These include all manifolds and in some sense are the stacks on which we can hope to do some geometry.  Further, certain morphisms between differentiable stacks, called \emph{representable morphisms}, are singled out as the ones to which we can ascribe familiar properties such as being surjective, a submersion, \emph{et cetera}.  It is natural to ask how the tangent stack functor affects differentiable stacks and representable morphisms.  \S\ref{DifferentiableSubsection} will recall the definition of differentiable stacks and representable morphisms in detail and will show that the tangent stack functor sends differentiable stacks to differentiable stacks and representable morphisms to representable morphisms.

Differentiable stacks can be represented by Lie groupoids, and every Lie groupoid represents a differentiable stack.  This allows one to give the following ad-hoc definition of the tangent stack of a differentiable stack \cite[4.6]{\Heinloth}.  Represent $\X$ by a Lie groupoid $\Gamma$ with structure-maps
\begin{equation}\label{StructureMaps}\xymatrix{
\Gamma_1\times_{\Gamma_0}\Gamma_1\ar[r]^-\mu & \Gamma_1\ar[r]^i & \Gamma_1\ar@<0.5ex>[r]^{s,t}\ar@<-0.5ex>[r] & \Gamma_0\ar[r]^e& \Gamma_1.
}\end{equation}
Take tangent bundles and derivatives everywhere to obtain a new Lie groupoid $T^\gpd\Gamma$ with spaces $T\Gamma_0$, $T\Gamma_1$ and structure maps
\begin{equation}\label{TangentStructureMaps}\xymatrix{
T\Gamma_1\times_{T\Gamma_0}T\Gamma_1\ar[r]^-{T\mu} & T\Gamma_1\ar[r]^{Ti} & T\Gamma_1\ar@<0.5ex>[r]^{Ts,Tt}\ar@<-0.5ex>[r] & T\Gamma_0\ar[r]^{Te}& T\Gamma_1
}\end{equation}
and then take the tangent stack of $\X$ to be the stack represented by $T^\gpd\Gamma$.  In \S\ref{TangentGroupoidsSubsection} we show that $T^\st\X$ is indeed the stack obtained from this construction.

Finally, in \S\ref{TangentColimitSubsection} we will describe $T^\st\X$ in terms of a colimit.  This may help the category-minded reader to visualize the tangent stack, and it is also an important component in proving some of the later results on the structure of tangent stacks.

In subsequent sections we will refer to $T^\st$ and $\pi^\st$ as simply $T\colon\StDiff\to\StDiff$ and $\pi\colon T\Rightarrow\Id_\StDiff$ respectively.

\subsection{Construction of the tangent stack functor.}\label{TangentConstructionSubsection}

In this section we construct the tangent stack functor.  This construction is just a stacky version of the construction of a geometric morphism between categories of sheaves from a morphism of sites.  See \cite[VII.10]{\MoerdijkMacLane}.

In what follows we will use arrows of the form $\to$, $\Rightarrow$, $\Rrightarrow$ to denote lax functors, lax natural transformations and modifications respectively.  The $2$-category of pseudofunctors $\mathcal{B}\to\mathcal{C}$ together with pseudonatural transformations and modifications will be denoted $[\mathcal{B},\mathcal{C}]$.  (The $2$-morphisms in all the categories we consider will be invertible, so that `lax' and `pseudo-' have the same meaning for us.)  $\Gpd$ denotes the $2$-category of groupoids.  See \cite[Chapter 7]{\Borceux} for the language of $2$-categories.

Let $i\colon\StDiff\hookrightarrow[\Diff^\mathrm{op},\Gpd]$ denote the inclusion of the $2$-category of stacks on $\Diff$ into the $2$-category of presheaves of groupoids on $\Diff$.  Precomposition with $T\colon\Diff\to\Diff$ determines a lax functor $T^\ast\colon [\Diff^\mathrm{op},\Gpd]\to [\Diff^\mathrm{op},\Gpd]$.

\begin{lemma}
$T^\ast$ restricts to a lax functor $T^\ast\colon \StDiff\to\StDiff$.
\end{lemma}
\begin{proof}
We must check that $T^\ast\colon[\Diff^\mathrm{op},\Gpd]\to[\Diff^\mathrm{op},\Gpd]$ sends stacks to stacks.  But $T\colon\Diff\to\Diff$ preserves open covers and pullbacks by open maps.  The stack condition for $T^\ast\Y$ now follows as an instance of the stack condition for $\Y$.
\end{proof}

Lax functors $F\colon\mathcal{C}\to\mathcal{D}$ and $G\colon\mathcal{D}\to\mathcal{C}$ are \emph{adjoint} ($F$ is left-adjoint to $G$, and $G$ is right-adjoint to $F$) if there is an equivalence of categories $\mathrm{Mor}_\mathcal{D}(Fc,d)\simeq \mathrm{Mor}_\mathcal{C}(c,Gd)$ lax natural in $c$ and $d$.  By this we mean that $\mathrm{Mor}_\mathcal{D}(F-,-)$ and $\mathrm{Mor}_\mathcal{C}(-,G-)$ are equivalent objects of $[\mathcal{C}^\mathrm{op}\times\mathcal{D},\Cat]$.

\begin{proposition}
$T^\ast\colon\StDiff\to\StDiff$ admits a left adjoint $T^\st\colon\StDiff\to\StDiff$ called the \emph{tangent stack functor}.
\end{proposition}

Left adjoints are determined up to natural equivalence, so the proposition defines the tangent stack functor.  Why should this left-adjoint be the functor we seek?  The functor $T^\ast$ is effectively determined by the equations
\[\mathrm{Mor}(X,T^\ast\Y)=\mathrm{Mor}(TX,\Y).\]
The fact that $T^\st$ is left-adjoint to $T^\ast$, however,  states that there is an equivalence
\[\mathrm{Mor}(\X,T^\ast\Y)\simeq\mathrm{Mor}(T^\st\X,\Y)\]
for any stack $\X$.  Thus $T^\st$ is determined by a property that, when restricted to manifolds, determines the tangent functor $T\colon\Diff\to\Diff$.  Everything else in this subsection will be a formal consequence of the adjunction of $T^\st$ with $T^\ast$.

\begin{proof}
We may assume that $\Diff$ is small.  Indeed, every object of $\Diff$ is isomorphic to a manifold embedded in some $\R^n$, so that $\Diff$ is equivalent to the full subcategory of $\Diff$ whose objects are these smooth manifolds embedded in some $\mathbb{R}^n$.  Moreover the $2$-category $\Gpd$ is cocomplete.  We may therefore form a left adjoint $T^\mathrm{pre}\colon[\Diff^\mathrm{op},\Gpd]\to[\Diff^\mathrm{op},\Gpd]$ to $T^\ast\colon [\Diff^\mathrm{op},\Gpd]\to [\Diff^\mathrm{op},\Gpd]$ by taking a left Kan extension.  There is also a left adjoint $a\colon[\Diff^\mathrm{op},\Gpd]\to\StDiff$ to $i\colon\StDiff\to[\Diff^\mathrm{op},\Gpd]$ given by sending a prestack to its associated stack.  Now $a\circ T^\mathrm{pre}\circ i$ is the required left adjoint:
\begin{eqnarray*}
\Mor(a\circ T^\mathrm{pre}\circ i\X,\Y)
&\simeq& \Mor(T^\mathrm{pre}\circ i\X,i\Y)\\
&\simeq&\Mor(i\X,T^\ast(i\Y))\\
&=&\Mor(i\X,i(T^\ast\Y))\\
&=&\Mor(\X,T^\ast\Y).
\end{eqnarray*}
This completes the proof.
\end{proof}

Now we wish to extend the the natural transformation $\pi\colon T\Rightarrow\Id_\Diff$, which consists of the projections $\pi_X\colon TX\to X$, to a lax natural transformation $\pi^\st\colon T^\st\Rightarrow\Id_\StDiff$. We will use the fact that for $2$-categories $\mathcal{B}$ and $\mathcal{C}$ the functor $[\mathcal{B},\mathcal{C}]\to[\mathcal{B}^\mathrm{op}\times\mathcal{C},\Cat]$, $F\mapsto\Mor_\mathcal{D}(F-,-)$ is locally full and faithful.  

\begin{definition}
Precomposition with $\pi\colon T\Rightarrow \Id_\Diff$ determines a natural transformation $\pi^\ast\colon\Id_\StDiff\Rightarrow T^\ast$.  The \emph{projection map} $\pi^\st\colon T^\st\Rightarrow\Id_\StDiff$ is the natural transformation corresponding to the composite
\[\Mor(\X,\Y)\xrightarrow{\pi^\ast\circ -}\Mor(\X,T^\ast\Y)\simeq\Mor(T^\st\X,\Y).\]
This means that there is a $2$-cell
\[\xymatrix{
\mathrm{Mor}(-,-)\ar[rr]^{\pi^\ast\circ -}_{}="1"\ar[rd]_{-\circ\pi^\st} & {} & \mathrm{Mor}(-,T^\ast-)\\
{} & \mathrm{Mor}(T^\st-,-)\ar[ur]_\simeq^{}="2"\ar@{=>}"2";"1" & {}
}\]
in $[\StDiff^\mathrm{op}\times\StDiff,\Cat]$.
\end{definition}

\subsection{Tangent stacks and the Yoneda embedding.}\label{YonedaSubsection}

In this subsection we show that, when restricted to manifolds using the Yoneda embedding, the tangent stack functor simply becomes the tangent functor and the natural projection $\pi^\st\colon T^\st\Rightarrow\Id$ becomes the projection $\pi\colon T\Rightarrow\Id$.

\begin{proposition}
There is a natural equivalence $\varepsilon\colon T^\st\circ y\Rightarrow y\circ T$.
\end{proposition}
\begin{proof}
There is an equivalence
\[\Mor(y(TX),\Y)=\Mor(yX,T^\ast\Y)\simeq\Mor(T^\st(yX),\Y)\]
natural in both variables.  Here the equality is the definition of $T^\ast\Y$ and the equivalence is from the adjunction of $T^\st$ with $T^\ast$.  But given $2$-categories $\mathcal{B}$ and $\mathcal{C}$, the functor $[\mathcal{B},\mathcal{C}]\to[\mathcal{B}^\mathrm{op}\times\mathcal{C},\Cat]$, $F\mapsto\Mor(F-,-)$ is locally full and faithful.  We therefore obtain the natural equivalence of the statement and a $2$-cell
\[\xymatrix{
\Mor(yT-,-)\ar@{=}[r]\ar@/^6ex/[rr]^{-\circ\varepsilon}^-{}="1" & \Mor(y-,T^\ast-)\ar@{=>}"1"\ar[r]^{\simeq}&\Mor(T^\st y-,-)
}\]
in $[\StDiff^\mathrm{op}\times\StDiff,\Cat]$.
\end{proof}

\begin{corollary}
Without loss of generality, we may assume that $T^\st\circ y = y\circ T$, which is to say that when restricted to $\Diff$, the tangent stack functor $T^\st\colon\StDiff\to\StDiff$ is just given by the tangent functor $T\colon\Diff\to\Diff$.
\end{corollary}

\begin{proposition}\label{NaturalTransformationProposition}
The two natural transformations $\pi^\st\ast\Id_y\colon T^\st\circ y\Rightarrow y$ and $\Id_y\ast\pi\colon y\circ T\Rightarrow y$ coincide under the identification  $T^\st\circ y=y\circ T$.  This means that the triangles
\[\xymatrix{
T^\st(yX)\ar[rd]_{\pi^\st_{yX}}\ar@{=}[rr] &{}& y(TX)\ar[dl]^{y(\pi_X)}\\
{}&yX &{}
}\]
commute, or even more simply, that when restricted to $\Diff$, $\pi^\st$ is given by $\pi$.
\end{proposition}
\begin{proof}
There are no nontrivial $2$-morphisms between morphisms between objects in the image of $y$.  Consequently, to show that the two natural transformations coincide it will suffice to show that there is a modification
\[\xymatrix{
T^\st\circ y\ar@{=>}[rd]_{\pi^\st\ast\Id_y}^{\phantom{.}}="1"\ar@{=}[rr] &{}& y\circ T\ar@{=>}[dl]^{\Id_y\ast\pi}_{\phantom{.}}="2"\\
{}&y.\ar@3{->}"1";"2" &{}
}\]
The lax natural transformations in this triangle determine a triangle
\[\xymatrix{
\Mor(T^\st\circ y-,-)\ar@{<-}[rd]_{-\circ(\pi^\st\ast\Id_y)}^{\phantom{.}}="1"\ar@{=}[rr] &{}& \Mor(y\circ T-,-)\ar@{<-}[dl]^{(\Id_y\ast\pi)\circ-}_{\phantom{.}}="2"\\
{}&\Mor(y-,-)&{}
}\]
in $[\StDiff^\mathrm{op}\times\StDiff,\Cat]$, and to construct the required modification it will suffice to fill this triangle with a $2$-cell.  But the triangle can be decomposed as three triangles
\[\xymatrix{
\Mor(T^\st\circ y-,-)\ar@/^6ex/@{=}[rr]_{}="1"\ar@{<-}[rd]_{-\circ(\pi^\st\ast\Id_y)}\ar@{<-}[r]^-\simeq &\Mor(y-,T^\ast-)\ar@{=}[r]& \Mor(y\circ T-,-)\ar@{<-}[dl]^{(\Id_y\ast\pi)\circ-}\\
{}&\Mor(y-,-)\ar[u]^{\pi^\ast\circ-}&{}
}\]
each of which can be filled with a $2$-cell.  The top triangle is filled with the $2$-cell obtained in the construction of $\varepsilon$ (which is assumed equal to  the identity), the left-hand triangle by the $2$-cell that defines $\pi^\st$, and the right-hand triangle commutes on the nose by definition.
\end{proof}

\subsection{Tangent stacks and tangent groupoids.}\label{TangentGroupoidsSubsection}

In this section we will prove that if a stack $\X$ is represented by a Lie groupoid $\Gamma$ then $T^\st\X$ is represented by the tangent Lie groupoid $T^\gpd\Gamma$.  In fact, we shall prove a much more precise functorial statement.

Let $\LieGpd$ denote the strict $2$-category of Lie groupoids and write
\[\B\colon\LieGpd\to\StDiff\]
for the lax functor that sends a Lie groupoid to its stack of torsors.  For a recollection on Lie groupoids we recommend \cite{\Moerdijk}, and for a definition of the stack $\B\Gamma$ of $\Gamma$-torsors we refer the reader to \cite[\S 2.4]{\BX}.  The promotion of the assignment $\Gamma\mapsto\B\Gamma$ to a lax functor can be written down directly, and is given by the composition of lax functors $\LieGpd\hookrightarrow\mathsf{Bi}$, $B\colon\mathsf{Bi}\to\StDiff$ described in \cite[\S 4]{\Lerman}.

\begin{definition}[Tangent groupoid functor]\label{TangentGroupoidDefinition}
Let  $T^\gpd\colon\LieGpd\to\LieGpd$ denote the strict functor that:
\begin{enumerate}
\item Sends a Lie groupoid $\Gamma$ with structure maps \eqref{StructureMaps} to the Lie groupoid $T^\gpd\Gamma$ with structure maps \eqref{TangentStructureMaps}.
\item Sends a morphism $f\colon\Gamma\to\Delta$ determined by maps $f_i\colon\Gamma_i\to\Delta_i$ to the morphism $T^\gpd f$ determined by the $Tf_i\colon T\Gamma_i\to T\Delta_i$.
\item Sends a $2$-morphism $\phi\colon f\Rightarrow g$ determined by $\phi\colon\Gamma_0\to\Delta_1$ to the $2$-morphism $T^\gpd\phi\colon T^\gpd f\Rightarrow T^\gpd g$ determined by $T\phi\colon T\Gamma_0\to T\Delta_1$.
\end{enumerate}
There is an evident natural morphism $\pi^\gpd\colon T^\gpd\Rightarrow\Id_\LieGpd$ obtained from the projection maps $T\Gamma_i\to\Gamma_i$.
\end{definition}

\begin{theorem}\label{TBCommuteTheorem}
There is a natural equivalence $T^\st\circ\B\simeq\B\circ T^\gpd$, which is to say, there are equivalences
\[T^\st(\B\Gamma)\simeq\B(T^\gpd\Gamma)\]
natural in $\Gamma$.  This equivalence identifies $\pi^\st_{\B\Gamma}$ with $\B\pi^\gpd_\Gamma$ in the sense that there is a modification
\[\xymatrix{
T^\st\circ\B \ar@{<=>}[rr]^{\simeq}\ar@{=>}[rd]_{\pi^\st}^{\phantom{.}}="1" & {} & \B\circ T^\gpd\ar@{=>}[dl]^{\B\pi^\gpd}_{\phantom{.}}="2"\\
{} & \B & {}\ar@3{->}"1";"2"
}\]
\end{theorem}
\begin{proof}
For a Lie groupoid $\Gamma$ and a stack $\Y$ let $\mathrm{Desc}(\Gamma,\Y)$ denote the groupoid whose objects are pairs $(f,\phi)$ consisting of a morphism $f\colon\Gamma_0\to\Y$ and a $2$-morphism $\phi\colon s^\ast f\Rightarrow t^\ast f$ for which $\pi_{23}^\ast \phi\circ\pi_{12}^\ast\phi=\pi_{13}^\ast\phi$ and whose arrows $\lambda\colon (f,\phi)\to(g,\psi)$ are $2$-morphisms $\lambda\colon f\Rightarrow g$ for which $\psi=t^\ast\lambda\circ\phi\circ s^\ast\lambda^{-1}$.

Then the stack condition and the fact that $\Gamma_0\to\B\Gamma$ is an atlas state that the $2$-commutative square 
\[\xymatrix{
\Gamma_1\ar[r]^s\ar[d]_t &\Gamma_0\ar[d]_{}="1"\\
\Gamma_0\ar[r]^{}="2" & \B\Gamma\ar@{=>}"1";"2"
}\]
determines an equivalence $\Mor(\B\Gamma,\Y)\to\mathrm{Desc}(\Gamma,\Y)$ \cite[2.20]{\BX};
this equivalence is natural in both variables.  Note that $\mathrm{Desc}(\Gamma,T^\ast\Y)=\mathrm{Desc}(T^\gpd\Gamma,\Y)$.  We therefore have an equivalence
\begin{eqnarray*}
\Mor(T^\st\B\Gamma,\Y)
&\simeq& \Mor(\B\Gamma,T^\ast\Y)\\
&\simeq&\mathrm{Desc}(\Gamma,T^\ast\Y)\\
&\cong&\mathrm{Desc}(T^\gpd\Gamma,\Y)\\
&\simeq&\Mor(\B T^\gpd\Gamma,\Y)
\end{eqnarray*}
natural in both variables.  The first result follows.  The second result can now be proved by carefully examining the sequence of equivalences above, just the modification was constructed in the proof of Proposition~\ref{NaturalTransformationProposition}.
\end{proof}

\subsection{Tangent stacks and differentiable stacks.}\label{DifferentiableSubsection}

We now recall the notions of differentiable stack and representable morphism and prove that the properties of differentiability and representability are preserved by the tangent stack functor $T^\st$.  The following definitions can be found in \cite{\BX} or \cite{\Heinloth}.

\begin{itemize}
\item A \emph{representable submersion} is a morphism $X\to\X$ with domain a manifold for which:
For any manifold $Y$ and any morphism $Y\to\X$, the fibre product $X\times_\X Y$ is representable and $X\times_\X Y\to Y$ is a submersion.
It is a \emph{representable {surjective} submersion} if in addition $X\times_\X Y\to Y$ is surjective.  Representable surjective submersions are also called \emph{atlases}.
\item A \emph{differentiable stack} is a stack on $\Diff$ that admits an {atlas}.
\item A morphism $\X\to\Y$ is \emph{representable} if:
For any representable submersion $Y\to\Y$, or for a single atlas $Y\to\Y$, the pullback $\X\times_\Y Y$ is representable.
It is called \emph{submersive, \'etale, proper} if, in addition, $\X\times_\Y Y\to Y$ is submersive, \'etale, proper.
\item An atlas $U\to\X$ yields a Lie groupoid $U\times_\X U\rightrightarrows U$ and an equivalence $\X\simeq\B(U\times_\X U\rightrightarrows U)$.  We say that $\X$ is \emph{represented by $U\times_\X U\rightrightarrows U$}.
\end{itemize}

\begin{theorem}\label{TangentStackProperties}
$T^\st$ sends differentiable stacks, representable morphisms, and representable (surjective) submersions to differentiable stacks, representable morphisms, and (surjective) submersions respectively.  If
\begin{equation}\label{CartesianSquare}\xymatrix{
\mathfrak{W}\ar[r]\ar[d]&\X\ar[d]_{}="1"\\
\Y\ar[r]^{}="2"&\mathfrak{Z}\ar@{=>}"1";"2"
}\end{equation}
is a cartesian diagram of differentiable stacks in which the morphisms are representable and one of $\Y\to\mathfrak{Z}$, $\X\to\mathfrak{Z}$ is a submersion, then the diagram
\[\xymatrix{
T^\st\mathfrak{W}\ar[r]\ar[d]&T^\st\X\ar[d]_{}="1"\\
T^\st\Y\ar[r]^{}="2"&T^\st\mathfrak{Z}\ar@{=>}"1";"2"
}\]
obtained by applying the lax functor $T^\st$ to \eqref{CartesianSquare} is again cartesian.
\end{theorem}
\begin{proof}
Let $\X$ be a differentiable stack.  Then $\X\simeq\B X$ for some groupoid $X$, and consequently $T^\st\X\simeq T^\st\B X\simeq \B T^\gpd X$, the second equivalence by Theorem~\ref{TBCommuteTheorem}.  Thus $T^\st\X$ is itself differentiable.

To show that $T^\st$ sends representable morphisms to representable morphisms we will use the fact that a morphism of Lie groupoids $f\colon\Gamma\to\Delta$ induces a representable morphism $\B f$ if and only if the map
\begin{gather*}
\label{RepresentabilityCriterion}\Delta_1\times_{\Delta_0}\Gamma_1\to(\Delta_1\times_{\Delta_0}\Gamma_0)\times(\Delta_1\times_{\Delta_0}\Gamma_0)\\
\nonumber(\delta,\gamma)\mapsto(\delta,s(\gamma))\times(\delta\cdot f_1(\gamma),t(\gamma))
\end{gather*}
is an embedding.

So let $f\colon\X\to\Y$ be representable.  By choosing an atlas for $\Y$ and taking the induced atlas for $\X$ we may find a diagram
\[\xymatrix{
\X\ar[r]^f\ar[d]_\simeq &\Y\ar[d]^\simeq_{}="1"\\
\B X\ar[r]_{\B f'}^{}="2" & \B Y\ar@{=>}"1";"2"
}\]
where $f'$ is a groupoid morphism satisfying the representability criterion above.  From Theorem~\ref{TBCommuteTheorem} we obtain a diagram
\[\xymatrix{
T^\st\X\ar[r]^{Tf}\ar[d] &T^\st\Y\ar[d]_{}="1"\\
T\B X\ar[r]_{T\B f'}^{}="2"\ar[d]_{} & T\B Y\ar[d]^{}_{}="3"\ar@{=>}"1";"2"\\
\B T^\gpd X\ar[r]_{\B T^\gpd f'}^{}="4" & \B T^\gpd Y\ar@{=>}"3";"4"
}\]
whose vertical maps are all equivalences, so that it will suffice to show that $\B T^\gpd f'$ is representable.  Since the map
\begin{gather*}
Y_1\times_{Y_0}X_1\to(Y_1\times_{Y_0}X_0)\times(Y_1\times_{Y_0}X_0)\\
(y,x)\mapsto(y,s(x))\times(y\cdot f'_1(x),t(x))
\end{gather*}
is an embedding and $T$ preserves pullbacks and embeddings, the map
\begin{gather*}
TY_1\times_{TY_0}TX_1\to(TY_1\times_{TY_0}TX_0)\times(TY_1\times_{TY_0}TX_0)\\
(y,x)\mapsto(y,s(x))\times(y\cdot Tf'_1(x),t(x))
\end{gather*}
is also an embedding.  It follows that $T^\gpd f'$ is representable, as required.

If $f$ is in addition a (surjective) submersion, then the component $f'_0\colon X\to Y$ could also be chosen a surjective submersion, so that $Tf'_0\colon TX\to TY$ is itself a (surjective) submersion, and then $\B T^\gpd f'$ is a (surjective) submersion also.

Finally consider the cartesian diagram \eqref{CartesianSquare}.  Choose a groupoid $Z$ representing $\mathfrak{Z}$, and then construct groupoids $X$, $Y$ representing $\X$ and $\Y$ by taking pullbacks.  We can form the pullback groupoid $X\times_Z W$, and $\B X\times_{\B Z}\B Y\simeq \B(X\times_Z W)$, so that $X\times_Z W$ represents $\mathfrak{W}$.  That is, the diagram \eqref{CartesianSquare} above is equivalent to one obtained by applying $\B$ to the cartesian diagram
\begin{equation}\label{GroupoidCartesianSquare}\xymatrix{
X\times_ZY\ar[r]\ar[d]&X\ar[d]_{}="1"\\
Y\ar[r]^{}="2"&{Z}\ar@{=>}"1";"2"
}\end{equation}
in $\LieGpd$.  Thus, by applying $T^\st$ to \eqref{CartesianSquare} we obtain a diagram that by Theorem~\ref{TBCommuteTheorem} is equivalent to applying $\B\circ T^\gpd$ to the diagram \eqref{GroupoidCartesianSquare}.  But it is simple to check that $T^\gpd(X\times_Z Y)=T^\gpd X\times_{T^\gpd Z}T^\gpd Y$, so that the diagram obtained by applying $\B\circ T^\gpd$ to \eqref{GroupoidCartesianSquare} is itself cartesian, as required.
\end{proof}

\subsection{Tangent stacks as lax colimits.}\label{TangentColimitSubsection}
In this last subsection we will show how to describe the tangent stack of a stack on $\Diff$ as a lax colimit.  This gives us a direct definition of $T^\st\X$ for any stack $\X$ on $\Diff$, regardless of whether $\X$ is differentiable, and gives us a description that is independent of a representing groupoid in that case.  See \cite[Chapter 7]{\Borceux} or \cite[Appendix 2]{\GH} for the definition of lax colimits.

Let $\X$ be a stack on $\Diff$. The \emph{category of manifolds over $\X$} is defined to be the comma category $(\Diff\downarrow\X)$.  An object in $(\Diff\downarrow\X)$ is simply a morphism
\begin{equation}\label{DiffXMorphism}W\to\X\end{equation}
whose domain is a manifold, and an arrow in $(\Diff\downarrow\X)$ from $W\to\X$ to $V\to\X$ is just a triangle
\begin{equation}\label{DiffXTriangle}\xymatrix@R=3pt{
W\ar[rrd]_(0.3){}="1"\ar[dd] &{}&{} \\
{}&{}& \X.\\
V\ar[rru]^(0.3){}="2"\ar@{=>}"1";"2"
}\end{equation}
Composition is given by pasting of diagrams.  There is an obvious strict functor
\[F_\X\colon(\Diff\downarrow\X)\to\StDiff\]
which remembers the manifolds in \eqref{DiffXMorphism} and \eqref{DiffXTriangle} but forgets the morphisms to $\X$.  There is also a tautological cone
\[c_\X\colon F_\X\Rightarrow\Delta_\X\]
determined by the morphisms in \eqref{DiffXMorphism} and the 2-morphisms in \eqref{DiffXTriangle}.

\begin{lemma}\label{StackColimitLemma}
The cone $c_\X$ determines an identification
\[\X=\colim F_\X\]
that we write informally as
\[\X=\colim_{W\to\X}W.\]
\end{lemma}
\begin{proof}
Since $\StDiff$ is a full subcategory of the functor category $[\Diff^\mathrm{op},\Gpd]$, the fact that composition with $c_\X$ determines an equivalence
$\mathrm{Mor}(\X,\Y)\xrightarrow{\simeq}\cone(F_\X,\Y)$ is an immediate consequence of the definitions.
\end{proof}

\begin{corollary}
Let $\X$ be a stack on $\Diff$.  Then
\[T^\st\X=\colim T\circ F_\X\]
or, informally
\[T^\st\X=\colim_{W\to\X}TW.\]
\end{corollary}
\begin{proof}
Since $T^\st$ is left adjoint to the functor $T^\ast$, it preserves colimits, and so $T^\st\X=T^\st\colim F_\X=\colim T^\st\circ F_\X=\colim T\circ F_\X$.  (We have suppressed the Yoneda embedding $y\colon\Diff\to\StDiff$ from our notation.)
\end{proof}

We will see in the sequel that this way of expressing the tangent stack can be very useful, since it gives us a way to describe the tangent stack $T^\st\X$ in terms of tangent bundles of manifolds without first having to choose a Lie groupoid representing $\X$.

\section{Vector Fields}\label{VectorFieldsSection}

This section extends the notion of vector field from manifolds to stacks on $\Diff$.  The definition is given in \S\ref{StackVFSubsection}.  We then show in \S\ref{VFSubmersionSubsection} that vector fields on stacks can be lifted through submersions; this is a technical result whose importance cannot be over-emphasised since it relates vector fields on a stack to vector fields on an atlas for that stack.  Then in \S\ref{GpdVFSubsection} we define vector fields on a Lie groupoid and show that they are equivalent to vector fields on the stack of torsors.  Finally \S\ref{VFSupportSubsection} defines the support of a vector field.  
\begin{definition}
A differentiable stack $\X$ is \emph{proper} if the diagonal $\Delta\colon\X\to\X\times\X$ is proper.  (The diagonal is always representable.)
\end{definition}
Some of the results in this section, and most of the results in the next section, are only proved for proper differentiable stacks.  Any manifold is a proper stack, as is any quotient by a compact Lie group.  Properness is best thought of as some sort of general Hausdorff or separability condition.  Appendix~\ref{ProperAppendix} recalls some properties of proper stacks in detail.  

In this section we will refer to the tangent stack functor and the projection map as $T\colon\StDiff\to\StDiff$ and $\pi\colon T\Rightarrow\Id_\StDiff$ respectively, rather than using the more elaborate notation of \S\ref{TangentStacksSection}.

\subsection{Vector fields on stacks.}\label{StackVFSubsection}

A vector field on a manifold $M$ is a section of the tangent bundle $TM$.  This means that a vector field is a map $X\colon M\to TM$ with the property that
\[\pi_M\circ X=\Id_M.\]
We wish to generalize this and define vector fields on any  stack on $\Diff$.  The ingredients are in place: any $\X$ has a tangent stack $T\X$ and a projection map $\pi_\X\colon T\X\to\X$.  However, we must bear in mind that within the $2$-category $\StDiff$ two morphisms can fail to be equal and yet still be isomorphic.  Indeed, the collection of morphisms $\X\to\Y$ is often vast when compared to its set of isomorphism classes, so to require that two morphisms be equal is quite unreasonable.  In particular, defining a vector field on $\X$ to be a morphism $\X\to T\X$ for which $\pi_\X\circ X=\Id_\X$ holds on the nose would not result in a useful notion.  Instead we weaken the equation to a $2$-morphism and define vector fields on stacks as follows.

\begin{definition}\label{VectorFieldDefinition}
Let $\X$ be a stack on $\Diff$.  A \emph{vector field on $\X$} is a pair $(X,a_X)$ consisting of a morphism
\[X\colon\X\to T\X\]
and a $2$-morphism $a_X\colon\pi_\X\circ X\Rightarrow\Id_\X$ that we depict in the diagram
\[\xymatrix{
\X\ar[r]^X\ar@/_6ex/[rr]_{\Id_\X}^{}="1"& T\X\ar[r]^{\pi_\X}\ar@{=>}"1"^{a_X}& \X.
}\]
\end{definition}

It is clear that if $\X$ is a manifold $M$ (or more correctly, the image of $M$ under the Yoneda embedding) then the vector fields on $\X$ form a set that is isomorphic to the set of vector fields on $M$.  However, the same comment that motivated the last definition --- that morphisms between stacks are very rarely equal but can still be isomorphic --- indicates that we should introduce a notion of isomorphism between vector fields on stacks, otherwise we may find ourselves dealing with an unmanageably large collection of vector fields.  Indeed, if $\X$ is equivalent to a manifold $M$ but not isomorphic to it, then the vector fields on $\X$ could form a collection far larger than the set of vector fields on $M$.  Our solution is the following.

\begin{definition}\label{EquivalenceDefinition}
Vector fields $X$ and $Y$ are \emph{equivalent} if there is $\lambda\colon X\Rightarrow Y$ for which $a_X=a_Y\circ(\Id_{\pi_\X}\ast\lambda)$.  We depict this relation as
\[\xymatrix{
\X\ar[r]^X\ar@/_6ex/[rr]^{}="1"& T\X\ar[r]\ar@{=>}"1"^{a_X}& \X&=&
\X\ar@/^3ex/[r]_{}="3"^>>>X\ar@/_3ex/[r]^{}="4"_>>>Y\ar@/_6ex/[rr]^{}="1"& T\X\ar[r]\ar@{=>}"1"^{a_Y}& \X.\ar@{=>}"3";"4"_{\lambda}
}\]
Such a $\lambda$ is called an \emph{equivalence}.  Vector fields and equivalences between them form the \emph{groupoid of vector fields on $\X$}, denoted $\Vect(\X)$.  We will often omit the anchoring $2$-morphisms $a_X$ from the notation, referring simply to vector fields $X$ on $\X$.
\end{definition}

With this definition one does find that the groupoid of vector fields on a representable stack $\X\simeq M$ is equivalent to the set of vector fields on $M$.  Indeed, we will see in Theorem~\ref{VectorFieldTheorem} that the groupoid of vector fields on a differentiable stack can be described easily in terms of vector fields on a Lie groupoid representing that stack.

\begin{example}[Manifolds]
If $\X=M$ is a manifold then a vector field on $\X$ is just a pair $(X,\Id)$ where $X$ is a vector field on $M$.   There are no nontrivial equivalences among these vector fields.  Thus $\Vect(\X)$ is just the set of vector fields on $M$.
\end{example}

\begin{example}[The zero vector field]
Recall that we can regard $\X$ as the colimit $\colim_{W\to\X}W$, and that $T\X$ is defined to be the colimit $\colim_{W\to\X}TW$.  The zero sections $W\to TW$ assemble into a natural transformation that induces a morphism
\[Z\colon\X\to T\X\]
on the colimits.  The fact that each composition $W\to TW\to W$ is the identity $\Id_W$ means that there is a uniquely-determined $2$-morphism
\[\xymatrix{
\X\ar[r]^Z\ar@/_6ex/[rr]_{\Id_\X}^{}="1"& T\X\ar[r]^{\pi_\X}\ar@{=>}"1"^{a_Z}& \X.
}\]
The pair $(Z,a_Z)$ is called the \emph{zero vector field} on $\X$.
\end{example}

\begin{example}[Products]
One consequence of Theorem~\ref{TangentStackProperties} is that, just as with manifolds, the derivatives of the projections $\X\times\Y\to\X$, $\X\times\Y\to\Y$ induce an equivalence
\begin{equation}\label{ProductEquivalence}T(\X\times\Y)\xrightarrow{\ \simeq\ }T\X\times T\Y.\end{equation}
By its construction this equivalence is compatible with the projections $\pi_{\X\times\Y}$ and $\pi_\X\times\pi_\Y$.  Consequently, if $(X,a_X)$, $(Y,a_Y)$ are vector fields on $\X$ and $\Y$ respectively, then we obtain a vector field $(X\times Y,a_X\times a_Y)$ on $\X\times\Y$.  (We are required to pick a quasi-inverse to the equivalence \eqref{ProductEquivalence}.)
\end{example}

\begin{example}[Differentiation with respect to time.]
The last example gives us an equivalence
\[T(\X\times\R)\cong T\X\times T\R\]
and in particular a vector field $\bigddt$ on $\X\times\R$ given by taking the product of the zero vector field on $\X$ and the unit vector field on $\R$.  This works just as well if $\R$ is replaced with an open interval $I\subset\R$.
\end{example}

\subsection{Vector fields and submersions.}\label{VFSubmersionSubsection}

The following lemma can be proved by a simple argument that uses partitions of unity and the fact that submersions of manifolds are locally projections.

\begin{lemma}\label{ManifoldLiftLemma}
Let $f\colon M\to N$ be a submersion of manifolds and let $X_N$ be a vector field on $N$.  Then there is a vector field $X_M$ on $M$ with the property that $Tf\circ X_M = X_N\circ f$:
\[\xymatrix{
M\ar[r]^{X_M}\ar[d]_{f} & TM\ar[d]^{Tf}\\
N\ar[r]_{X_N}&TN
}\]
\end{lemma}

This subsection will extend the lemma above from submersions of manifolds to representable submersions of differentiable stacks.  This is a necessary step if we are to get a handle on vector fields on stacks.  For the concrete way to understand a differentiable stack is to choose an atlas, thus representing the stack by a Lie groupoid.  But an atlas is just a representable surjective submersion, and so an appropriate generalization of Lemma~\ref{ManifoldLiftLemma} would allow us to take a vector field on a differentiable stack $\X$ and `lift' it to a vector field on an atlas $U$ for $\X$.  We make this generalization below and exploit it in \S\ref{GpdVFSubsection}, where the groupoid of vector fields on $\X$ is described explicitly in terms of a Lie groupoid representing $\X$.

\begin{lemma}\label{VectorFieldLemma}
Let $\Y$ be a proper differentiable stack, let $s\colon\Y\to\X$ be a representable submersion and let $(X_\X,a_\X)$ be a vector field on $\X$.  Then we may find a vector field $(X_\Y,a_\Y)$ on $\Y$ and a commutative diagram
\begin{equation}\label{VectorFieldLemmaSquare}\xymatrix{
\Y\ar[r]^{X_\Y}\ar[d] & T\Y\ar[d]_{}="1"\\
\X\ar[r]_{X_\X}^{}="2"&T\X\ar@{=>}"1";"2"
}\end{equation}
for which the $2$-morphisms in
\begin{equation}\label{VectorFieldLemmaConditionSquare}\xymatrix{
\Y\ar@/^6ex/[rr]_{}="6"\ar[r]\ar[d] & T\Y\ar@{=>}"6"_{a_\Y}\ar[d]_{}="1"\ar[r] &\Y\ar[d]_{}="3"\\
\X\ar@/_6ex/[rr]^{}="5" \ar[r]^{}="2"&T\X\ar@{=>}"5"^{a_\X}\ar@{=>}"1";"2" \ar[r]^{}="4" & \X\ar@{=>}"3";"4"
}\end{equation}
compose to give the trivial $2$-morphism from $s\colon\Y\to\X$ to itself.  In the square on the right the horizontal maps are the projections $\pi_\X\colon T\X\to\X$, $\pi_\Y\colon T\Y\to\Y$, and the square itself is obtained from the lax naturality of the projection $\pi\colon T\Rightarrow\Id_\StDiff$.
\end{lemma}

This lemma is a direct generalization of Lemma~\ref{ManifoldLiftLemma}, for it reduces to that lemma in the case that $\X$ and $\Y$ are manifolds.  However, it contains a significant new feature in the condition on diagram~\eqref{VectorFieldLemmaConditionSquare}, which is vacuous in the manifold case.  This is a typical feature of categorification, but why does it arise?  One answer is to consider what the condition means: the traditional diagram \eqref{VectorFieldLemmaSquare} relates the {morphisms} $X_\X$ and $X_\Y$, but the new diagram \eqref{VectorFieldLemmaConditionSquare} relates the {vector fields} $(X_\X,a_\X)$ and $(X_\Y,a_\Y)$.  A better answer, of course, is that this condition is useful.  It means that $X_\X$ and $X_\Y$ induce a new vector field on the pullback $\Y\times_\X\Y$, in the following sense.

\begin{lemma}\label{EssentialSurjectivityLemma}
Suppose that we are in the situation of Lemma~\ref{VectorFieldLemma}.  The diagram \eqref{VectorFieldLemmaSquare} induces a morphism
\[X_{\Y\times_\X \Y}\colon \Y\times_\X \Y\longrightarrow T\Y\times_{T\X}T\Y=T(\Y\times_\X \Y).\]
that is a vector field on $\Y\times_\X\Y$.  In other words, there is a $2$-morphism $\pi_{\Y\times_\X\Y}\circ X_{\Y\times_\X\Y}\Rightarrow\Id_{\Y\times_\X\Y}$.
\end{lemma}

For the purposes of later reference we record Lemma~\ref{VectorFieldLemma} in the special case that $\Y$ is a manifold.  This corollary is essential for our applications.

\begin{corollary}\label{VectorFieldCorollary}
Let $U\to\X$ be a representable submersion and let $X$ be a vector field on $\X$. Then we may find a vector field $X_U$ on $U$ and a commutative diagram
\begin{equation}\label{VectorFieldSquare}\xymatrix{
U\ar[r]^{X_U}\ar[d] & TU\ar[d]_{}="1"\\
\X\ar[r]_X^{}="2"&T\X\ar@{=>}"1";"2"
}\end{equation}
for which the $2$-morphisms in
\begin{equation}\label{VectorFieldConditionSquare}\xymatrix{
U\ar[r]\ar[d] & TU\ar[d]_{}="1"\ar[r] &U\ar[d]_{}="3"\\
\X\ar@/_6ex/[rr]^{}="5" \ar[r]^{}="2"&T\X\ar@{=>}"5"^{a_X}\ar@{=>}"1";"2" \ar[r]^{}="4" & \X\ar@{=>}"3";"4"
}\end{equation}
compose to give the identity.
\end{corollary}

\begin{proof}[Proof of Lemma~\ref{VectorFieldLemma}.]
It is possible to construct a $2$-commutative diagram
\begin{equation}\label{VectorFieldProofDiagram}\xymatrix{
\mathfrak{W}\ar[r]\ar[d] & T\Y\ar[dr]\ar[d]_{}="1" & {} \\
\Y\ar[r]^{}="2"\ar[d] & \Y\times_\X T\X\ar[r]\ar[d]_{}="3" & \Y\ar[d]_{}="4"\\
\X\ar[r]^{}="5" & T\X\ar[r]^{}="6" & \X \ar@{=>}"1";"2"\ar@{=>}"4";"6"\ar@{=}"3";"5"
}\end{equation}
as follows.  The bottom-right square is obtained by taking pullbacks.
The morphism $T\Y\to\Y\times_\X T\X$ is determined by the square
\[\xymatrix{
T\Y\ar[r]^{\pi_\Y}\ar[d] & \Y\ar[d]_{}="1"\\
T\X\ar[r]_{\pi_\X}^{}="2" & \X,\ar@{=>}"1";"2"}\]
which is to say, from naturality of $\pi$.  The morphism $\Y\to \Y\times_\X T\X$ is determined by the square
\[\xymatrix{
\Y\ar@{=}[r]\ar[d]_{X_{\X}|\Y} & \Y\ar[d]_{}="1"\\
T\X\ar[r]^{}="2" &\X.\ar@{=>}"2";"1"^{a_\X|\Y}
}\]
Strict commutativity of the triangle is now immediate, as is strict commutativity of the bottom-left square.  The composition in the middle row is just $\Id_\Y$, and if the $2$-morphisms in the bottom two squares are pasted with $a_\X$ we recover the trivial $2$-morphism.  The top-left square is simply a pullback, with $\mathfrak{W}$ shorthand for $T\Y\times_{\Y\times_\X T\X}\Y$.

We claim that we can find a morphism $s$ and a $2$-morphism $\sigma$ in a diagram of the form
\begin{equation}\label{WeakSectionDiagram}\xymatrix{
\Y\ar@/_6ex/[rr]^{}="5"_{\Id_\Y} \ar[r]^{s}="2"&\mathfrak{W}\ar[r]\ar@{=>}"5"^{\sigma}& \Y.
}\end{equation}
Assuming this for the time being, the vector field $(X_\Y,a_\Y)$ can now be constructed using $s$, $\sigma$, and the top half of diagram \eqref{VectorFieldProofDiagram}.  The required diagram \eqref{VectorFieldLemmaSquare} can be constructed using the weak section $s$ and the left-hand part of \eqref{VectorFieldProofDiagram}.  The composition of the $2$-morphisms in the resulting diagram \eqref{VectorFieldLemmaConditionSquare} can now be computed directly and seen to be trivial.  This proves the lemma.

We now show how to construct diagram \eqref{WeakSectionDiagram}.  To do this we will first study the morphism $T\Y\to\Y\times_\X T\X$, from which $\mathfrak{W}\to\Y$ is obtained by pulling back.  There are equivalences
\begin{eqnarray*}
T\Y&\simeq& T\Y\times_{T\X}T\X\\
{}&=& T\Y\times_{T\X}\colim TW\\
{}&\simeq&\colim (T\Y\times_{T\X}TW)\\
{}&\simeq&\colim T(\Y\times_\X W)\\
\Y\times_\X T\X &=& \Y\times_\X\colim TW\\
{}&\simeq&\colim(\Y\times_\X TW)
\end{eqnarray*}
Here the lax colimits are all taken over the category $(\Diff\downarrow\X)$ of morphisms $W\to\X$ with $W$ a manifold.  We have used the fact that lax colimits in $\StDiff$ commute with pullbacks and that $T$ preserves pullbacks under submersions (Theorem~\ref{TangentStackProperties}).  One can check from the chains of equivalences given above that there is a $2$-morphism in the square
\[\xymatrix{
T\Y\ar[r]^-{\simeq}\ar[d] & \colim T(\Y\times_\X W)\ar[d]_{\phantom{.}}="1"\\
\Y\times_\X T\X\ar[r]_-{\simeq}^-{\phantom{.}}="2"&\colim \Y\times_\X TW,\ar@{=>}"1";"2"
}\]
or in other words that $T\Y\to \Y\times_\X T\X$ is the colimit of the projections $T(\Y\times_\X W)\to \Y\times_\X TW$.

Since $\Y\to\X$ is a representable submersion, each $\Y\times_\X W$ is a manifold and $\Y\times_\X W\to W$ is a submersion.  Thus each of the projections $T(\Y\times_\X W)\to\Y\times_\X TW$ is a fibrewise-linear surjection of vector bundles over the manifold $\Y\times_\X W$.  In particular, each of these projections is in a natural way an \emph{affine vector bundle}, i.e.~a fibre bundle with fibres isomorphic to $\R^n$ and with structure group $\R^n\rtimes\mathrm{GL}(n,\R)$, where $\R^n$ acts on itself by translation.  In this case
\begin{eqnarray*}
n&=&\mathrm{dim}(\Y\times_\X W)-\mathrm{dim}(W)\\
&=&\mathrm{dim}(\Y)-\mathrm{dim}(\X).
\end{eqnarray*}
Thus each $T(\Y\times_\X W)\to \Y\times_\X TW$ is an affine vector bundle.  What is more, in the diagram
\[\xymatrix{
T(\Y\times_\X W_1)\ar[d]\ar[r]& T(\Y\times_\X W_2)\ar[d]\\
\Y\times_\X TW_1\ar[r]&W\times_\X TW_2
}\]
induced by a morphism in $(\Diff\downarrow\X)$ the horizontal maps constitute a morphism of affine vector bundles.  The morphism $T\Y\to \Y\times_\X T\X$ is then the colimit of a diagram of affine vector-bundles, and so is itself an affine vector-bundle whose base is a stack.  This is a simple consequence of the fact that colimits in stacks commute with pullbacks.  Finally, since $\mathfrak{W}\to\Y$ is obtained from this morphism by pulling back, it is itself an affine vector-bundle, this time with base the proper stack $\Y$.

We have shown that $\mathfrak{W}\to\Y$ is an affine vector-bundle and we wish to show that it admits a weak section, i.e.~construct diagram \eqref{WeakSectionDiagram}.  Locally, $\Y$ has the form $[M/G]$ with $G$ compact, and an affine vector-bundle on $[M/G]$ is a $G$-equivariant affine vector bundle on $M$ (Corollary \ref{LocalQuotientCorollary}).  This bundle on $M$ admits a section, and averaging with respect to $G$ we can assume that the section is $G$-invariant, so that the bundle on $[M/G]$ admits a section.  So locally we can find sections of an affine vector-bundle on $\Y$.  These can be glued using a partition of unity (Proposition \ref{POUProposition}) to obtain the required global section.
\end{proof}

\begin{proof}[Proof of Lemma~\ref{EssentialSurjectivityLemma}]
The composition $\pi_{\Y\times_\X \Y}\circ X_{\Y\times_\X \Y}$  is the morphism $\Y\times_\X \Y\to \Y\times_{\X}\Y$ determined by the following square:
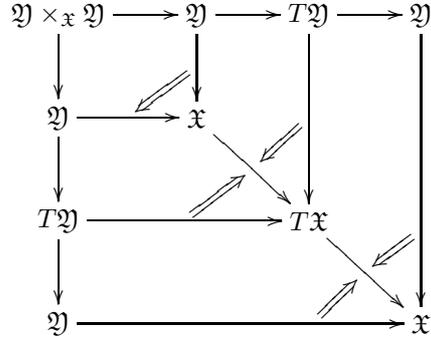
\[\xymatrix{
\Y\times_{\X}\Y\ar[r]\ar[d] & \Y\ar[d]_{}="1"\ar[r] & T\Y\ar[dd]_{}="3"\ar[r] & \Y\ar[ddd]_(0.7){}="5"\\
\Y\ar[r]^{}="2"\ar[d]& \X\ar[rd]^{}="9"_{}="10" &{} & {}\\
T\Y\ar[rr]^{}="4"\ar[d] & {} & T\X\ar[rd]^{}="6"_{}="7" & {}\\
\Y\ar[rrr]^(0.7){}="8" & {}& {}& \X
\ar@{=>}"1";"2"\ar@{=>}"3";"9"\ar@{=>}"4";"10"\ar@{=>}"5";"6"\ar@{=>}"8";"7"
}\]
We can paste copies of $a_\Y$ onto the upper and left edges to obtain a new square
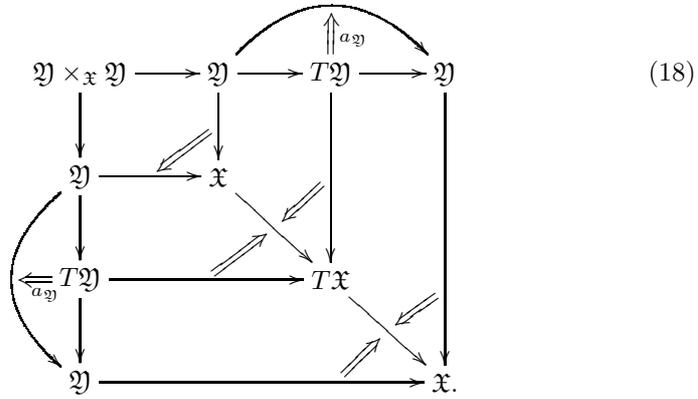
\begin{equation}\label{BigSquareTwo}\xymatrix{
\Y\times_{\X}\Y\ar[r]\ar[d] & \Y\ar[d]_{}="1"\ar@/^6ex/[rr]_{}="11"\ar[r] & T\Y\ar@{=>}"11"_{a_\Y}\ar[dd]_{}="3"\ar[r] & \Y\ar[ddd]_(0.7){}="5"\\
\Y\ar@/_6ex/[dd]^{}="12"\ar[r]^{}="2"\ar[d]& \X\ar[rd]^{}="9"_{}="10" &{} & {}\\
T\Y\ar@{=>}"12"^{a_\Y}\ar[rr]^{}="4"\ar[d] & {} & T\X\ar[rd]^{}="6"_{}="7" & {}\\
\Y\ar[rrr]^(0.7){}="8" & {}& {}& \X.
\ar@{=>}"1";"2"\ar@{=>}"3";"9"\ar@{=>}"4";"10"\ar@{=>}"5";"6"\ar@{=>}"8";"7"
}\end{equation}
The effect of this modification is to replace the morphism $\Y\times_\X\Y\to\Y\times_\X\Y$ by a new morphism that is related to the original by a $2$-morphism determined by $a_\Y$.  To prove the result it therefore remains to show that the morphism $\Y\times_\X\Y\to\Y\times_\X\Y$ determined by the new square \eqref{BigSquareTwo} is the identity.  We can perform a simple manipulation on the square --- inserting a copy of $a_\X$ directly adjacent to a copy of its inverse --- without altering the composite $2$-morphism, to obtain a new square:
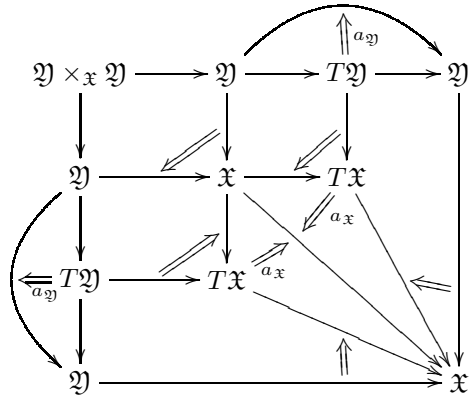
\[\xymatrix{
\Y\times_{\X}\Y\ar[r]\ar[d] & \Y\ar[d]_{}="1"\ar@/^6ex/[rr]_{}="13"\ar[r] & T\Y\ar@{=>}"13"_{a_\Y}\ar[d]_{}="3"\ar[r] & \Y\ar[ddd]_(0.7){}="5"\\
\Y\ar@/_6ex/[dd]^{}="14"\ar[r]^{}="2"\ar[d]& \X\ar[d]_{}="10" \ar[r]^{}="9"\ar[rrdd]^(0.3){}="11"_(0.3){}="12" & T\X \ar[rdd]^{}="6"\ar@{=>}"11"^{a_\X}& {}\\
T\Y\ar@{=>}"14"^{a_\Y}\ar[r]^{}="4"\ar[d] &  T\X\ar[rrd]_{}="7"\ar@{=>}"12"_{a_\X}& {}& {}\\
\Y\ar[rrr]^(0.7){}="8" &{} &{} & \X
\ar@{=>}"1";"2"\ar@{=>}"3";"9"\ar@{=>}"4";"10"\ar@{=>}"5";"6"\ar@{=>}"8";"7"
}\]
This new square contains two copies of the diagram \eqref{VectorFieldLemmaConditionSquare}.  The condition of Lemma~\ref{VectorFieldLemma} now means that we may replace each copy of \eqref{VectorFieldLemmaConditionSquare} with the much simpler
\[\xymatrix{
\Y\ar@{=}[r]\ar[d] & \Y\ar[d]\\
\X\ar@{=}[r] & \X.
}\]
Our modified version of \eqref{BigSquareTwo} now simplifies to give the standard pullback square so that the morphism $\Y\times_\X \Y\to \Y\times_\X \Y$ is just the identity, as required.
\end{proof}

\subsection{Vector fields on Lie groupoids.}\label{GpdVFSubsection}

In Definition~\ref{TangentGroupoidDefinition} we defined the tangent groupoid functor $T^\gpd\colon\LieGpd\to\LieGpd$ that assigns to each Lie groupoid its tangent groupoid.  We will write $T^\gpd$ as $T$ for simplicity.

\begin{definition}[Vector fields on a Lie groupoid]\label{LieGpdVectorFieldDefinition}
Let $\Gamma$ be a Lie groupoid.   A \emph{vector field on $\Gamma$} is a groupoid morphism $X\colon\Gamma\to T\Gamma$ for which the composition $\pi_\Gamma\circ X$ is the identity on $\Gamma$.  An \emph{equivalence} between vector fields $X$, $Y$ on $\Gamma$ is a $2$-morphism $\psi\colon X\Rightarrow Y$ for which $\Id_{\pi_\Gamma}\ast\psi=\Id_{\Id_\Gamma}$.  The vector fields on $\Gamma$ and equivalences between them together define the \emph{groupoid of vector fields on $\Gamma$}, denoted $\Vect(\Gamma)$.
\end{definition}

Note the relative simplicity of the definition of vector fields on a groupoid in comparison with that of vector fields on stacks, Definition~\ref{VectorFieldDefinition}.  We have asked that the composite $\pi_\Gamma\circ X$ be {equal to} the identity on $\Gamma$, not that it is merely $2$-isomorphic to the identity.   This would have been the wrong choice for stacks since morphisms are so rarely equal, but for groupoids it is the correct notion, as we see in the following theorem.

\begin{theorem}\label{VectorFieldTheorem}
The groupoid of vector fields on a Lie groupoid $\Gamma$ is equivalent to the groupoid of vector fields on the stack $\B\Gamma$.
\end{theorem}

It may appear that this theorem is a simple consequence of the Dictionary Lemma~\ref{DictionaryLemma}, which tells us how to relate stack morphisms $\B\Gamma\to T\B\Gamma\simeq\B T\Gamma$ to Lie groupoid morphisms $\Gamma\to T\Gamma$.  However, the Dictionary Lemma only tells us about those morphisms $\B\Gamma\to\B T\Gamma$ that we already know can be lifted to a morphism $\Gamma_0\to T\Gamma_0$.  The essential ingredient, then, is Corollary~\ref{VectorFieldCorollary}, which guarantees that any vector field on $\B\Gamma$ does admit such a lift.

\begin{proof}[Proof of Theorem~\ref{VectorFieldTheorem}.]
This result is proved by combining functoriality of $\B$ and $T$ to construct a functor $\Vect(\Gamma)\to\Vect(\B\Gamma)$ and then using Corollary~\ref{VectorFieldCorollary} and the Dictionary Lemma~\ref{DictionaryLemma} in order to prove that the functor is an equivalence.  

Lax functoriality of $\B$ together with the lax natural equivalence and modification of Theorem~\ref{TBCommuteTheorem} determine, for each vector field $X$ on $\Gamma$, $2$-cells
\[\xymatrix{
\B\Gamma\ar[rr]^{\B X}\ar[dr]_{\B\Id_\Gamma}^{}="2" &{}& \B T\Gamma\ar[dl]^{\B\pi_\Gamma}_{}="1"\\
&\B\Gamma\ar@{=>}"1";"2"&
}\]
\[\xymatrix{
\B\Gamma\ar@/^2ex/[rr]^{\B\Id_\Gamma}_{}="1"\ar@/_2ex/[rr]_{\Id_{\B\Gamma}}^{}="2" &{}& \B\Gamma\ar@{=>}"1";"2"
}
\qquad
\xymatrix{
\B T\Gamma\ar[rr]^{\simeq}\ar[rd]^{}="1"_{\B\pi_\Gamma} &{}& T\B\Gamma\ar[dl]_{}="2"^{\pi_{\B\Gamma}}\\ 
&\B\Gamma &\ar@{=>}"1";"2"
}\]
These three diagrams may be pasted together to obtain a new diagram
\[\xymatrix{
\B\Gamma\ar[r]^{BX}\ar@/_6ex/[rr]_{\Id_{\B\Gamma}}^{}="1"& T\B\Gamma\ar[r]^{\pi_{\B\Gamma}}\ar@{=>}"1"^{b_X}& \B\Gamma.
}\]
We therefore have an assignment $X\mapsto (BX,b_X)$ from vector fields on $\Gamma$ to vector fields on $\B\Gamma$.  Standard properties of lax functors, lax natural transformations and modifications allow us to promote this assignment to a functor $\Vect(\Gamma)\to\Vect(\B\Gamma)$.  We will prove the theorem by showing that this functor is an equivalence.

It is possible to verify from the construction above that we have commutative diagrams
\begin{equation}\label{LieGpdVFProofDiagramOne}
\xymatrix{\Gamma_0\ar[r]^{X_0}\ar[d] & T\Gamma_0\ar[d]_{}="1"\\
\B\Gamma\ar[r]_{{B}X}^{}="2"&T\B\Gamma\ar@{=>}"1";"2"}
\end{equation}
satisfying the conclusion of Corollary~\ref{VectorFieldCorollary} and with the further property that the induced map $\Gamma_1\to T\Gamma_1$ --- which by Lemma~\ref{EssentialSurjectivityLemma} is itself a vector field ---- is just $X_1$.

Part \ref{DictionaryThree} of the Dictionary Lemma~\ref{DictionaryLemma}, combined with diagram \eqref{LieGpdVFProofDiagramOne}, immediately show that equivalences $BX\Rightarrow BY$ are in correspondence with equivalences $X\Rightarrow Y$.  Thus the functor $\Vect(\Gamma)\to\Vect(\B\Gamma)$ is fully faithful.  We now wish to show that it is essentially surjective.  Let $\tilde X$ be a vector field on $\B\Gamma$.  Since $\Gamma_0\to\B\Gamma$ is an atlas we may apply Corollary~\ref{VectorFieldCorollary} to obtain a vector field $X_0$ on $\Gamma_0$ and a diagram
\begin{equation}\label{LieGpdVFProofDiagramTwo}
\xymatrix{\Gamma_0\ar[r]^{X_0}\ar[d] & T\Gamma_0\ar[d]_{}="1"\\
\B\Gamma\ar[r]_-{\tilde X}^{}="2"&T\B\Gamma\ar@{=>}"1";"2"}
\end{equation}
satisfying the conclusion of Corollary~\ref{VectorFieldCorollary}.  Then, by Lemma~\ref{EssentialSurjectivityLemma}, the induced map $X_1\colon\Gamma_1\to T\Gamma_1$ is itself a vector field and the pair $X_0$, $X_1$ together define a vector field on $\Gamma$.  But now we may compare diagrams \eqref{LieGpdVFProofDiagramOne} and  \eqref{LieGpdVFProofDiagramTwo} and, using part \ref{DictionaryTwo} of Dictionary Lemma~\ref{DictionaryLemma}, conclude that there is a $2$-morphism $\tilde X\Rightarrow BX$.  Since the diagrams \eqref{LieGpdVFProofDiagramOne} and  \eqref{LieGpdVFProofDiagramTwo} satisfied the conclusion of Corollary~\ref{VectorFieldCorollary} it follows that $\tilde X\Rightarrow BX$ is in fact an equivalence of vector fields.  Thus the functor $\Vect(\Gamma)\to\Vect(\B\Gamma)$ is essentially surjective.
\end{proof}

\subsection{The support of a vector field.}\label{VFSupportSubsection}

The \emph{support} of a vector field on a manifold is the closure of the set of points on which the vector field is nonzero.  Equivalently, the support of a vector field on a manifold is the complement of the largest open set on which the vector field vanishes.  Extending this notion to stacks presents a problem: vector fields are very rarely equal to zero (by which we mean, equal to the zero vector field) but might more often be \emph{equivalent} to zero.  We therefore wish to consider the `largest open substack on which the vector field is equivalent to the zero vector field'.  In order to do so we will prove in this subsection that such a largest open substack exists, provided that the stack admits smooth partitions of unity, as is the case with all proper stacks (see Definition~\ref{POUDefinition} and Proposition~\ref{POUProposition}).

\begin{proposition}\label{PatchEquivalenceProposition}
Let $X$ and $Y$ be vector fields on a differentiable stack $\X$ that admits partitions of unity.  If $X$ and $Y$ are equivalent on full open substacks $\mathfrak{A}_\alpha$ of $\X$, then they are equivalent on the full open substack $\bigcup\mathfrak{A}_\alpha$.  In particular, there is a unique maximal open substack of $\X$ on which $X$ and $Y$ are equivalent.
\end{proposition}

\begin{definition}
Let $X$ be a vector field on a differentiable stack $\X$ that admits partitions of unity.  Then the \emph{support of $X$} is defined to be the subset $\mathrm{supp}(X)$ of $\bar\X$ whose complement corresponds to the largest open substack of $\X$ on which $X$ is equivalent to the zero vector field.
\end{definition}

Here $\bar\X$ refers to the \emph{underlying space} of $\X$, which is canonically homeomorphic to the orbit space of any Lie groupoid representing $\X$.  Open subsets of $\bar\X$ correspond to full open substacks of $\X$ \cite[\S 2]{\Morse}.

\begin{proof}[Proof of Proposition~\ref{PatchEquivalenceProposition}.]
Let $U\to\X$ be an atlas and take vector fields $X_U$, $Y_U$ on $U$ and commutative diagrams
\[\xymatrix{
U\ar[r]^{X_U}\ar[d] & TU\ar[d]_{}="1"\\
\X\ar[r]_{X}^{}="2"& T\X\ar@{=>}"1";"2"
}\qquad
\xymatrix{
U\ar[r]^{Y_U}\ar[d] & TU\ar[d]_{}="1"\\
\X\ar[r]_{Y}^{}="2"& T\X\ar@{=>}"1";"2"
}\]
satisfying the conclusion of Corollary~\ref{VectorFieldCorollary}.  Given a full open substack $\mathfrak{A}$ of $\X$, we write $U_\mathfrak{A}$ for the subset of $U$ whose points lie in $\mathfrak{A}$.  Now by part~\ref{DictionaryThree} of the Dictionary Lemma~\ref{DictionaryLemma} we find that equivalences of vector fields $\lambda\colon X|\mathfrak{A}\Rightarrow Y|\mathfrak{A}$ are in 1-1 correspondence with maps
\[l\colon U_\mathfrak{A}\to T(U_\mathfrak{A}\times_\X U_\mathfrak{A})\]
which have the properties
\begin{enumerate}
\item $l$ is a lift of the unit map $U_\mathfrak{A}\to U_\mathfrak{A}\times_\X U_\mathfrak{A}$.
\item $T\pi_1\circ l = X_U|U_\mathfrak{A}$ and $T\pi_2\circ l=Y_U|U_\mathfrak{A}$.
\item The two composites
\[\xymatrix{U_{\mathfrak{A}}\times_\X U_{\mathfrak{A}}\ar[rr]^-{X_{U\times_\X U}\times_{TU}l} &{}& T(U_{\mathfrak{A}}\times_\X U_{\mathfrak{A}}\times_\X U_{\mathfrak{A}})\ar[r]^-{\pi_1\times_{\X}\pi_3} & T(U_{\mathfrak{A}}\times_\X U_{\mathfrak{A}})}\]
\[\xymatrix{U_{\mathfrak{A}}\times_\X U_{\mathfrak{A}}\ar[rr]^-{l\times_{TU}Y_{U\times_\X U}} &{}& T(U_{\mathfrak{A}}\times_\X U_{\mathfrak{A}}\times_\X U_{\mathfrak{A}})\ar[r]^-{\pi_1\times_{\X}\pi_3} & T(U_{\mathfrak{A}}\times_\X U_{\mathfrak{A}})}\]
coincide.
\end{enumerate}

The key to the proof is that a collection of maps $l$ satisfying the above conditions can be `patched' to obtain a new map that still satisfies the conditions; in other words, these conditions are preserved under averaging.

Define $\bar{\mathfrak{B}}\subset\bar{\X}$ to be the union $\bigcup\bar{\mathfrak{A}}_\alpha$.  This $\bar{\mathfrak{B}}$ is open, and we define $\mathfrak{B}$ to be the corresponding open substack of $\X$.  We will prove the proposition by constructing a $2$-morphism of vector fields $\lambda\colon X|\mathfrak{B}\Rightarrow Y|\mathfrak{B}$.

Since $\mathfrak{B}$ admits partitions of unity we may take a countable family of morphisms $\phi_i\colon\X\to\R$, with values in $[0,1]$, such that $\bar\phi_i\colon\bar\X\to\R$ is a partition of unity, and such that each $\phi_i$ is supported in one of the substacks $\mathfrak{A}_{\alpha_i}$.  Take an equivalence $\lambda_i\colon X|\mathfrak{A}_{\alpha_i}\Rightarrow Y|\mathfrak{A}_{\alpha_i}$ and write $l_i\colon U_{\mathfrak{A}_{\alpha_i}}\to T(U_{\mathfrak{A}_{\alpha_i}}\times_\X U_{\mathfrak{A}_{\alpha_i}})$ for the corresponding map.  Write $\phi_i l_i\colon U_{\mathfrak{A}_{\alpha_i}}\to T(U_{\mathfrak{A}_{\alpha_i}}\times_\X U_{\mathfrak{A}_{\alpha_i}})$ for the product of $l_i$ with the composition $U_{\mathfrak{A}_{\alpha_i}}\to\X\xrightarrow{\phi_i}\R$.  Since the supports of the $\bar\phi_i$ form a locally-finite family on $\bar\X$, the supports of the $\phi_i l_i$ also form a locally-finite family, and so we may form the sum $l=\sum \phi_i l_i\colon U_\mathfrak{B}\to T(U_\mathfrak{B}\times_\X U_\mathfrak{B})$.  It is now immediate to verify from its construction that $l$ satisfies conditions 1, 2 and 3 above, and so corresponds to an equivalence $\lambda\colon X|\mathfrak{B}\Rightarrow Y|\mathfrak{B}$ as required.
\end{proof}

\section{Integrals and Flows}\label{IntegralsFlowsSection}

Let $X$ be a vector field on a manifold $M$.  Recall that an \emph{integral curve} of $X$ through $m\in M$ is a curve $\gamma$ in $M$ such that $\gamma(0)=m$ and $\dot\gamma(t)=X(\gamma(t))$, which we can write as
\begin{equation}\label{IntegralCurveEquation}
T\gamma\circ\bigddt=X\circ\gamma.
\end{equation}
The \emph{flow} of $X$ is a smooth map $\phi\colon M\times\R\to M$
such that each $\phi(m,-)$ is the integral curve of $X$ through $m$.

\begin{proposition*}[{Existence and uniqueness of integral curves, \cite[I, 1.5]{\KobayashiNomizuOne}, \cite[2.4]{\Milnor}}]\label{ManifoldIntegralProposition}\hfill
\begin{enumerate}
\item The integral curve of $X$ through $m$ is unique where it is defined.
\item Integral curves exist for small time and depend smoothly on their initial value.  That is, for each $m_0\in M$ there is an open neighbourhood $U$ of $m_0$, an $\epsilon>0$, and a smooth map $\phi\colon U\times(-\epsilon,\epsilon)\to M$ such that each $\phi(m,-)$ is an integral curve of $X$ through $m$.
\end{enumerate}
\end{proposition*}

\begin{proposition*}[{Existence and uniqueness of flows, \cite[I, 1.6]{\KobayashiNomizuOne}}]\label{ManifoldFlowProposition}
If the flow of $X$ exists then it is unique.  If $X$ is compactly supported, then the flow of $X$ does exist.
\end{proposition*}

This section extends the notion of integral curve and flow from manifolds to stacks on $\Diff$, and proves analogues for proper differentiable stacks of the existence and uniqueness results above (they can fail if the manifold is not proper).  In all cases what one sees are weakened, or categorified, forms of the usual definitions and results, with equations replaced by $2$-morphisms, and with new conditions on, and relations among, these $2$-morphisms.  We begin in \S\ref{DefinitionsSubsection} with the definitions of integral morphisms and flows, then the existence and uniqueness results are stated in \S\ref{StatementsSubsection}, and finally the proofs are given in \S\ref{ProofsSubsection}.

\subsection{Definitions.}\label{DefinitionsSubsection}

Throughout this section $I$ will denote an open interval in $\R$.  We allow $I$ to be infinite, e.g.~$I=(0,\infty)$ or $I=\R$.

\begin{definition}[Integral morphisms]\label{IntegralDefinition}
Let $X$ be a vector field on a stack $\X$ on $\Diff$.  Then $\Phi\colon\Y\times I\to\X$ is \emph{an integral morphism of $X$} if there is a $2$-morphism
\begin{equation}\label{IntegralTwoMorphism}
t_\Phi\colon X\circ\Phi\Longrightarrow T\Phi\circ\bigddt,
\end{equation}
which we represent as the diagram
\begin{equation}\label{FlowDiagram}\xymatrix{
T(\Y\times I)\ar[r]^-{T\Phi}_{}="2"& T\X\\
\Y\times I\ar[u]^{\bigddt}\ar[r]_-\Phi &\X.\ar[u]_X^{}="1"\ar@{=>}^{t_\Phi}"1";"2"
}\end{equation}
The $2$-morphism $t_\Phi$ must satisfy the property that the $2$-morphisms in
\begin{equation}\label{CategorifiedFlowDiagram}\xymatrix{
\Y\times I\ar@{<-}@/_12ex/[dd]^{}="1"_{\Id_{\Y\times I}}\ar[r]_{}="3"^\Phi & \X\ar@{<-}@/^12ex/[dd]_{}="2"^{\Id_{\Y\times I}} \\
T(\Y\times I)\ar[u]\ar[r]^-{T\Phi}_{}="4"\ar@{=>}"1"^-{a_\smallddt} & T\X\ar@{=>}"2"_-{a_X}\ar[u]_{}="5"\\
\Y\times I\ar[u]\ar[r]_\Phi & \X\ar[u]_{}="6"\ar@{=>}"6";"4"^{t_\Phi}\ar@{=>}"5";"3"
}\end{equation}
compose to the trivial $2$-morphism from $\Phi\colon\Y\times\R\to\X$ to itself.  The choice of $t_\Phi$ is regarded as part of the data for $\Phi$.  Note that if $\Phi$ integrates $X$ and there is an equivalence $\lambda\colon Y\Rightarrow X$, then $\Phi$ also integrates $Y$ when equipped with the $2$-morphism $t_\Phi\circ (\lambda\ast\Id_\Phi)$.
\end{definition}

Consider the definition above when $\Y=\mathrm{pt}$ and $\X$ is a manifold.  The existence of $t_\Phi$ simply becomes the original equation \eqref{IntegralCurveEquation} while the condition on diagram \eqref{CategorifiedFlowDiagram} becomes vacuous. We therefore recover the definition of integral curves.  In general though, there may be many different choices of $t_\Phi$, only some of which satisfy the condition on diagram \eqref{CategorifiedFlowDiagram}.  This new condition, however, is a necessary one.  For we know from Corollary~\ref{VectorFieldCorollary} that a vector field $X$ on a stack $\X$ may be lifted to a vector field $X_U$ on an atlas $U$ for $\X$, and the new condition is what will allow us to relate the integral morphism $\Phi$ to the integral curves of $X_U$ on $U$.

In extending the uniqueness of integral curves to stacks we cannot expect that an integral morphism $\Phi\colon\Y\times I\to\X$ is determined by its initial value $\Phi|\Y\times\{0\}$, as is the case for integral curves on manifolds.  What we can ask is that the initial value determines $\Phi$ up to a $2$-morphism.  Indeed, this will be the case and the $2$-morphism in question will be uniquely determined, so long as we ensure that it satisfies the conditions in the next definition.

\begin{definition}[Integral $2$-morphisms]\label{IntegralTwoMorphismDefinition}
Let $\Phi,\Psi\colon\Y\times I\to\X$ integrate $X$.  An \emph{integral $2$-morphism} is a $2$-morphism $\Lambda\colon\Phi\Rightarrow\Psi$ that respects $t_\Phi$ and $t_\Psi$ in the sense that $(T\Lambda\ast\Id_{\smallddt})\circ t_\Phi = t_\psi\circ (\Id_{X}\ast\Lambda)$, which we express in diagrams as
\[\xymatrix{
T(\Y\times I)\ar[r]^{T\Phi}="4"_{}="2"\ar@/^8ex/[r]^{T\Psi}_{}="3"& T\X\\
\Y\times I\ar[u]\ar[r]_-\Phi &\X\ar[u]^{}="1"\ar@{=>}^{t_\Phi}"1";"2"
\ar@{=>}"4";"3"_{T\Lambda}}
\quad = \quad \xymatrix{
T(\Y\times I)\ar[r]^-{T\Psi}_{}="6"& T\X\\
\Y\times I\ar@/_8ex/[r]_{\Phi}^{}="8"\ar[u]\ar[r]_{\Psi}="7" &\X.\ar[u]^{}="5"\ar@{=>}^{t_\Psi}"5";"6"\ar@{=>}"8";"7"_{\Lambda}
}\]
\end{definition}

\begin{definition}[Flows]\label{FlowDefinition}
Let $X$ be a vector field on a stack $\X$ on $\Diff$.  A \emph{flow of $X$} is a morphism
\[\Phi\colon\X\times\R\to\X\]
integrating $X$ and equipped with a $2$-morphism $e_\Phi\colon\Phi|\X\times\{0\}\Rightarrow\Id_\X$.
\end{definition}

The isomorphism $e_\Phi$ in the last definition is simply our weakening of the initial condition on the integral curve through a point on a manifold. We will see that although flows are not unique, they are determined up to an integral $2$-morphism that is itself determined by $e_\Phi$.

\subsection{Existence and uniqueness theorems.}\label{StatementsSubsection}

\begin{theorem}[Uniqueness of integrals]\label{UniquenessTheorem}
Let $\X$ and $\Y$ be differentiable stacks and let $X$ be a vector field on $\X$.  Let $\Phi,\Psi\colon\Y\times I\to\X$ be morphisms that integrate $X$.  Then:
\begin{enumerate}
\item\label{UniquenessOne}
If $\Lambda,M\colon\Phi\Rightarrow\Psi$ are integral $2$-morphisms that coincide when restricted to some $\Y\times\{t_0\}$, then $\Lambda=M$.
\item\label{UniquenessTwo}
If $\X$ is proper, then any $2$-morphism $\lambda\colon\Phi|\Y\times\{t_0\}\Rightarrow\Psi|\Y\times\{t_0\}$ extends to a unique integral $2$-morphism $\Lambda\colon\Phi\Rightarrow\Psi$.
\end{enumerate}
\end{theorem}

This theorem is our generalization of the uniqueness of integral curves.  The next theorem is our generalization of the existence of integral curves.  We would like to say that the integral curve of $X$ through any point of $\X$ exists for small time, or more generally that any morphism $\phi\colon\Y\to\X$ extends to an integral morphism $\Phi\colon\Y\times (-\epsilon,\epsilon)\to\X$ that restricts to $\phi$ at time zero.  Of course, we must weaken this requirement slightly:

\begin{theorem}[Existence of integral morphisms]\label{ExistenceTheorem}
Let $X$ be a vector field on a proper differentiable stack $\X$.  Let $\Y$ be differentiable and let $\phi\colon\Y\to\X$ be a morphism whose image has compact closure.  Then for some $\epsilon>0$ there is a morphism $\Phi\colon\Y\times (-\epsilon,\epsilon)\to\X$ integrating $X$ and a $2$-morphism $\Phi|\Y\times\{0\}\Rightarrow\phi$.
\end{theorem}

\begin{note}
Both of Theorems~\ref{UniquenessTheorem} and \ref{ExistenceTheorem} can fail if one does not assume that $\X$ is proper.  Examples demonstrating this are given in \cite[\S 5.3]{\Morse}.
\end{note}

One is often interested in representable morphisms of differentiable stacks since these are, roughly speaking, the morphisms to which we can ascribe geometric properties.  The next result tells us when an integral morphism is representable.  The theorem is trivial when restricted to manifolds, since all maps of manifolds are representable.

\begin{theorem}[Representability of integrals]\label{RepresentabilityTheorem}
Let $X$ be a vector field on a proper differentiable stack $\X$, let $\Phi\colon\Y\times I\to\X$ integrate $X$, and fix any $t_0\in I$.  Then $\Phi$ is representable if and only if $\Phi|\Y\times\{t_0\}$ is representable.
\end{theorem}

Finally we extend the existence and uniqueness of flows to proper differentiable stacks.

\begin{theorem}[Existence and uniqueness of flows]\label{FlowTheorem}
Let $X$ be a vector field on a proper differentiable stack $\X$. 
\begin{enumerate}
\item A flow of $X$, if it exists, is unique up to a uniquely-determined integral $2$-morphism.  More precisely, if $\Phi$ and $\Psi$ are two flows of $X$, then there is a unique integral $2$-morphism $\Lambda\colon\Phi\Rightarrow\Psi$ such that $\Lambda|\X\times\{0\}=e_\Psi^{-1}e_\Phi$.
\item A flow of $X$, if it exists, is representable.
\item If $X$ has compact support, then a flow of $X$ does exist.
\end{enumerate}
\end{theorem}

\subsection{Proofs.}\label{ProofsSubsection}

\begin{definition}
Let us establish some notation.  Let $f\colon U\to V$ be a smooth map between manifolds and let $U,V$ carry vector fields $X_U$, $X_V$ respectively.  Then we say that $f$ \emph{intertwines $X_U$ and $X_V$}, or that \emph{$X_U$ and $X_V$ are compatible}, if $Tf\circ X_U = X_V\circ f$.
\end{definition}

\begin{proof}[Proof of Theorem~\ref{UniquenessTheorem}, part~\ref{UniquenessOne}]
First note that, if $Y\to\Y$ is an atlas, then $\Phi| Y\times I$, $\Psi|Y\times I$ still integrate $X$, and that $\Lambda|Y\times I$, $M|Y\times I$ are $2$-morphisms of these integrals that coincide on $Y\times\{t_0\}$.  Thus, if the conclusion of Theorem~\ref{UniquenessTheorem}  part~\ref{UniquenessOne} holds for manifolds, then $\Lambda|Y\times I=M|Y\times I$, and so $\Lambda=M$.  We may therefore assume that $\Y=Y$ is a manifold.  We may also without loss assume that $\Psi=\Phi$.

Let $U\to\X$ be an atlas and choose a vector field $X_U$ on $U$, with a diagram
\[\xymatrix{
U\ar[r]^{X_U}\ar[d]& TU\ar[d]_{}="1"\\
\X\ar[r]_X^{}="2" & T\X \ar@{=>}"1";"2"
}\]
satisfying the conclusion of Corollary~\ref{VectorFieldCorollary}.  Let $V\to Y\times I$ be the atlas obtained in the pullback-diagram
\begin{equation}\label{SetupDiagramTwo}\xymatrix{
V\ar[d]\ar[r]^{\tilde\Phi} & U\ar[d]_{}="1"\\
Y\times I\ar[r]_\Phi^{}="2" & \X. \ar@{=>}"1";"2"
}\end{equation}
By taking pullbacks in the rows of
\begin{equation}\label{UniquenessProofEquation}\xymatrix{
T(Y\times I)\ar[r]^-{T\Phi}_{}="1" & T\X & TU\ar[l]_{}="2"\\
Y\times I\ar[u]^{\bigddt}\ar[r]_-\Phi & \X\ar[u]^{}="3"_{}="4"_<<<X & U\ar[l]\ar[u]_{X_U}
\ar@{=>}"2";"3"\ar@{=>}"4";"1"^{t_\Phi}}\end{equation}
we obtain a smooth map $V\to TV$ that by the condition on diagram \eqref{CategorifiedFlowDiagram} and the condition of Corollary~\ref{VectorFieldCorollary} is itself a vector field $X_V$ on $V$; the proof is a mild generalization of the proof of Lemma~\ref{EssentialSurjectivityLemma}.
This vector field is compatible with $\smallddt$ via $V\to Y\times I$ and with $X_U$ via $\Phi\colon V\to U$.

After these preparations we can apply Dictionary Lemma~\ref{DictionaryLemma} part~\ref{DictionaryThree} to conclude that $\Lambda$ and $M$ determine and are determined by maps
\[l,m\colon V\to U\times_\X U.\]
These maps have the following properties:
\begin{enumerate}
\item Write $V_0$ for the part of $V$ that lies over $Y\times\{t_0\}$.  Then $l$ and $m$ coincide when restricted to $V_0$.  This is a consequence of the fact that $\Lambda$ and $M$ coincide when restricted to $Y\times\{t_0\}$.
\item $l$ and $m$ intertwine $X_V$ and $X_{U\times_\X U}$.  This follows from the  construction of $l$ and $m$ and the fact that $\Lambda$ and $M$ are integral $2$-morphisms; the proof involves diagram manipulations of the sort made in Lemma~\ref{EssentialSurjectivityLemma} and is left to the reader.
\item If $v_1,v_2\in V$ have equal images in $Y\times I$, then $l(v_1)=m(v_1)$ if and only if $l(v_2)=m(v_2)$.  This is because the assumption yields $\alpha\in U\times_\X U$ such that 
\begin{gather*}
l(v_1)\cdot\alpha = \alpha\cdot l(v_2),\\
m(v_1)\cdot\alpha=\alpha\cdot m(v_2).
\end{gather*}
Here $\cdot$ denotes composition in the groupoid $U\times_\X U\rightrightarrows U$.
\end{enumerate}

We now show that $l=m$.  Since $l$ and $m$ determine $\Lambda$ and $M$ respectively it will follow that $\Lambda=M$ as required.  Let $v\in V$, lying over $(v_Y,t_1)\in Y \times I$, and assume without loss that $t_1>t_0$.  Consider the map $I\to Y\times I$, $t\mapsto(v_Y,t)$, and the corresponding pullback diagram
\[\xymatrix{
\tilde I\ar[r]\ar[d] & V\ar[d]\\
I\ar[r]& Y\times I
}\]
There is a vector field on $\tilde I$ compatible with $\bigddt$ on $I$ with $X_V$ on $V$.  Since $\tilde I\to I$ is a surjective submersion we may therefore write $[t_0,t_1]=[s_0,s_1]\cup\cdots\cup[s_{n-1},s_n]$ and find $\gamma\colon[s_{i-1},s_i]\to V$ integrating $X_V$ and such that $\gamma_1(s_0)\in V_0$, such that $\gamma_i(s_i)$ and $\gamma_{i+1}(s_{i})$ have equal images in $Y\times I$, and such that $\gamma_n(s_n)$ and $v$ have equal images in $Y\times I$.  Now $l(\gamma_1(s_0))=m(\gamma_1(s_0))$ by the first property above.  Then $l(\gamma_1(s_1))=m(\gamma_1(s_1))$ by the second property above, so that $l(\gamma_2(s_1))=m(\gamma_2(s_1))$ by the third property above.  Continuing in this way we can conclude that $l(v)=m(v)$.  Since $v$ was chosen arbitrarily, $l=m$ as required.
\end{proof}

We now move onto the proof of Theorem~\ref{UniquenessTheorem} part~\ref{UniquenessTwo}.  This part requires us to construct the $2$-morphism and requires the additional condition of properness. It is consequently significantly more difficult than the proof of part~\ref{UniquenessOne}.  We begin by proving a series of lemmas and then assembling the proof from these.  The only aspect of properness that we use is the result of the first of these lemmas below.

\begin{lemma}\label{ProperLemma}
Let $M$, $N$ be smooth manifolds equipped with vector fields $X_M$, $X_N$ respectively.  Let $\pi\colon M\to N$ be a smooth proper map that intertwines $X_M$ and $X_N$.  Fix $m\in M$.  If the integral curve of $X_N$ through $\pi(m)$ exists to time $t$, then so does the integral curve of $X_M$ through $m$.
\end{lemma}
\begin{proof}
Suppose not.  Let $\gamma\colon[0,t]\to N$ denote the integral curve of $X_N$ with $\gamma(0)=\pi(m)$.  Without loss assume that the integral curve of $X_M$ through $m$ can be defined on $[0,t)$ but not on $[0,t]$.  Since $\pi$ is proper we may find some neighbourhood $U$ of $\gamma(t)$ and some $\epsilon>0$ such that the integral curve of $X_M$ through any point of $\pi^{-1}(U)$ can be defined on the interval $(-\epsilon,\epsilon)$.  So now choose $s\in(t-\epsilon,t)$ large enough that $\gamma(s)\in U$.  Then $\delta(s)\in\pi^{-1}(U)$, so that $\delta$ can be defined on $[0,s+\epsilon)$, which includes $t$.  This is a contradition.  This concludes the proof.
\end{proof}

Before we state the next lemma consider the following.  Suppose we are in the situation of Theorem~\ref{UniquenessTheorem}, part~\ref{UniquenessTwo}.  Let $p\colon\Y'\to\Y$ be a surjective submersion, and write also $p\colon\Y'\times I\to\Y\times I$ for the product of $p$ with the identity.  Then $\Phi\circ p,\Psi\circ p\colon\Y'\times I\to\X$ are both integrals of $X$, and the $2$-morphism $\Lambda$, if it existed, would induce an integral $2$-morphism $\Lambda'=\Lambda\ast\Id_{p}\colon\Phi\circ p\Rightarrow\Psi\circ p$ extending $\lambda\ast\Id_{p}$.  The converse is also true:

\begin{lemma}\label{UniquenessLemmaOne}
The conclusion of Theorem~\ref{UniquenessTheorem} part~\ref{UniquenessTwo} holds if there exists an integral $2$-morphism $\Lambda'\colon\Phi\circ p\Rightarrow\Psi\circ p$.
\end{lemma}
\begin{proof}
Let $P$ denote the $2$-morphism in the cartesian diagram
\[\xymatrix{
\Y'\times_\Y\Y'\times I\ar[r]^{\pi_1}\ar[d]_{\pi_2} &\Y'\times I\ar[d]^{p}_{}="1"\\
\Y'\times I\ar[r]^{}="2"_{p} & \Y\times I.\ar@{=>}"1";"2"
}\]
Then $\Lambda'$ descends to the required $\Lambda$ if the two composite $2$-morphisms
\[\xymatrix{
\Phi\circ p\circ\pi_1\ar@{=>}[r]^-{\pi_1^\ast(\Lambda')}&\Psi\circ p\circ\pi_1\ar@{=>}[r]^{\Psi_\ast P}&\Psi\circ p\circ\pi_2
}\]
\[\xymatrix{
\Phi\circ p\circ\pi_1\ar@{=>}[r]^-{\Phi_\ast P}&\Psi\circ p\circ\pi_2\ar@{=>}[r]^{\pi_2^\ast(\Lambda')}&\Psi\circ p\circ\pi_2
}\]
coincide.  But these are integral $2$-morphisms and, since $\Lambda'|\Y'\times\{t_0\}$ descends to $\lambda\circ\Id_p$, they coincide when restricted to $\Y'\times\{t_0\}$.  Then by Theorem~\ref{UniquenessTheorem} part~\ref{UniquenessOne} the composites coincide, as required.
\end{proof}

\begin{lemma}\label{UniquenessLemmaTwo}
Let $X$ be a vector field on a proper differentiable stack $\X$.  Let $Y$ be a manifold and let $\Phi,\Psi\colon Y\times I\to\X$ integrate $X$.  Then for any $Y_1\subset Y$ open with compact closure, and for any $t_1\in I$, we can find an open interval $J$ containing $t_1$ and contained in $I$, with the property that any $\lambda\colon\Phi|Y\times\{t_1\}\Rightarrow\Psi|Y_1\times\{t_1\}$ with $t_1\in J$ extends to an integral $2$-morphism $\Lambda\colon\Phi|Y_1\times J\Rightarrow\Psi|Y_1\times J$.
\end{lemma}
\begin{proof}
Let $U\to\X$ be an atlas and let $X_U$ be a vector field on $U$ satisfying the conclusions of Corollary~\ref{VectorFieldCorollary}.  Then, as in the proof of Theorem~\ref{UniquenessTheorem} part~\ref{UniquenessOne}, we can find an atlas $V\to Y\times I$ and a commutative diagram \eqref{SetupDiagramTwo} where $V$ carries a vector field $X_V$ compatible with $X_U$ on $U$ and $\smallddt$ on $Y\times I$.  Write $V_{t_1}$ for the part of $V$ that lies over $Y\times\{t_1\}$.  Around each point of $Y$ we can find an open neighbourhood small enough to lift to $V_{t_1}$ and small enough that the integral of $X_V$ through this lift exists on some small time interval $J$ containing $t_1$.  Since $Y_1$ has compact closure we may therefore find a cover $W_\Phi\to Y_1$ and a commutative diagram 
\[\xymatrix{
W_\Phi\times J\ar[r]^{\tilde\Phi}\ar[d] & U\ar[d]_{}="1"\\
Y_1\times J\ar[r]_\Phi^{}="2"&\X\ar@{=>}"1";"2"
}\]
in which $\tilde\Phi$ integrates $X_U$.  Repeating this process for $\Psi$, reducing $J$ and refining the two covers of $Y_1$ if necessary, we obtain a single cover $W\to Y_1$ and commutative diagrams
\[\xymatrix{
W\times J\ar[r]^{\Phi_1}\ar[d]& U\ar[d]_{}="1"\\
Y_1\times J\ar[r]_{\Phi|}^{}="2" &\X\ar@{=>}"1";"2"
}
\qquad
\xymatrix{
W\times J\ar[r]^{\Psi_1}\ar[d]& U\ar[d]_{}="1"\\
Y_1\times J\ar[r]_{\Psi|}^{}="2" &\X\ar@{=>}"1";"2"
}\]
in which $\Phi_1$, $\Psi_1$ both integrate $X_U$.  We no longer require the assumption on $\mathrm{cl} Y_1\subset Y$ and so by Lemma~\ref{UniquenessLemmaOne} we may replace $Y_1$ with $W$ and so assume that $\Phi|$ and $\Psi|$ factorize as follows:
\begin{equation}\label{UniquenessLemmaTwoEquationOne}\xymatrix{
& U\ar[d]_{}="1"\\
Y_1\times J\ar[r]_{\Phi|}^{}="2"\ar[ur]^{\Phi_1} &\X\ar@{=>}"1";"2"
}\qquad
\xymatrix{
& U\ar[d]_{}="1"\\
Y_1\times J\ar[r]_{\Psi|}^{}="2"\ar[ur]^{\Psi_1} &\X\ar@{=>}"1";"2"
}\end{equation}

The vector field $X_V$ on $V$ was constructed by pulling back in the rows of diagram \eqref{UniquenessProofEquation}.  It possible to use this fact to check from the construction of the diagrams \eqref{UniquenessLemmaTwoEquationOne} that composing the $2$-morphisms in the diagrams
\begin{equation}\label{UniquenessLemmaTwoEquationTwo}\xymatrix{
T(Y_1\times J)\ar[r]^-{T\Phi_1}\ar@/^6ex/[rr]_(0.57){}="1" & TU\ar[r]_{}="2"\ar@{=>}"1" & T\X\\
Y_1\times J\ar[u]^{\bigddt}\ar[r]^-{\Phi_1}\ar@/_6ex/[rr]^(0.57){}="3" & U\ar[u]_{X_U}\ar[r]\ar@{=>}"3" & \X\ar[u]_X^{}="4" \ar@{=>}"2";"4"}
\qquad
\xymatrix{
T(Y_1\times J)\ar[r]^-{T\Psi_1}\ar@/^6ex/[rr]_(0.57){}="1" & TU\ar[r]_{}="2"\ar@{=>}"1" & T\X\\
Y_1\times J\ar[u]^{\bigddt}\ar[r]^-{\Psi_1}\ar@/_6ex/[rr]^(0.57){}="3" & U\ar[u]_{X_U}\ar[r]\ar@{=>}"3" & \X\ar[u]_X^{}="4" \ar@{=>}"2";"4"
}\end{equation}
yields
\[\xymatrix{
T(Y_1\times J)\ar[r]^-{T\Phi|}_{}="1" & T\X\\
Y_1\times J\ar[u]^{\bigddt}\ar[r]_-{\Phi|}&\X\ar[u]_X^{}="2"\ar@{=>}"2";"1"^{t_\Phi}
}
\qquad
\xymatrix{
T(Y_1\times J)\ar[r]^-{T\Psi|}_{}="1" & T\X\\
Y_1\times J\ar[u]^{\bigddt}\ar[r]_-{\Psi|}&\X\ar[u]_X^{}="2"\ar@{=>}"2";"1"^{t_\Psi}
}\]
respectively.

With these preparations we will now prove that the conclusion of the lemma holds with the current choice of open interval $J$.  Let $\lambda\colon\Phi|Y\times\{t_1\}\Rightarrow\Psi|Y\times\{t_1\}$.  Using diagrams \eqref{UniquenessLemmaTwoEquationOne} and Dictionary Lemma~\ref{DictionaryLemma} part~\ref{DictionaryThree}, this $\lambda$ determines and is determined by a map $l\colon Y_1\times\{t_1\}\to U\times_\X U$.  By construction, this $l$ is a lift of $\Phi_1|\times\Psi_1|\colon Y_1\times\{t_1\}\to U\times U$.  By Lemma~\ref{ProperLemma}, since $\Phi_1\times\Psi_1$ integrates $X_U\times X_U$, and the proper map $\pi_1\times\pi_2\colon U\times_\X U\to U\times U$ intertwines $X_{U\times_\X U}$ and $X_U\times X_U$, we may form $L\colon Y_1\times J\to U\times_\X U$ integrating $X_{U\times_\X U}$ and restricting to $l$ on $Y_1\times\{t_1\}$.  By construction $\pi_1\circ L=\Phi_1$, $\pi_2\circ L=\Psi_1$, and so we obtain $\Lambda\colon\Phi|\to\Psi|$ extending $\lambda|$ by composing the $2$-morphisms in the diagram
\[\xymatrix{
{}&{}&U\ar[drr]_{}="1"&{}&{}\\
Y_1\times J\ar[r]^L\ar@/^12ex/[rrrr]^{\Phi|}_(0.485){}="3"\ar@/_12ex/[rrrr]_{\Psi|}^(0.485){}="4" & U\times_\X U\ar[dr]_{}="6"\ar[ur]^{}="5"&{} &{}& \X\\
{}&{}& U\ar[urr]^{}="2" & {}&{}\ar@{=>}"6";"4"\ar@{=>}"5";"3"\ar@{=>}"1";"2"
}\]

It remains to check that $\Lambda$ is an integral $2$-morphism, i.e.~that $(T\Lambda\ast\Id_{\smallddt})\circ t_\Phi = t_\psi\circ (\Id_{X}\ast\Lambda)$.  To do so we may use the description of $t_\Phi$, $t_\Psi$ given in diagram \eqref{UniquenessLemmaTwoEquationTwo} and the construction of $\Lambda$ to see --- after some tedious manipulation of diagrams as in the proof of Lemma~\ref{EssentialSurjectivityLemma} --- that the required result follows from the fact that $L$ integrates $X_{U\times_\X U}$.  
\end{proof}

\begin{proof}[Proof of Theorem~\ref{UniquenessTheorem}, part~\ref{UniquenessTwo}]
First, by Lemma~\ref{UniquenessLemmaOne} we may assume that $\Y=Y$ is a manifold.

Let $Y_1\subset Y$ be an open subset with compact closure.  We will construct an integral $2$-morphism $\Phi|Y_1\times I\Rightarrow\Psi|Y_1\times I$ extending $\lambda|Y_1$.  Given an open interval $J$ containing $t_0$ and contained in $I$ we will write $\Lambda_J\colon\Phi|Y_1\times J\Rightarrow\Psi|Y_1\times J$ for the unique integral $2$-morphism extending $\lambda|Y_1$, if it exists.  Applying Lemma~\ref{UniquenessLemmaTwo} with $t_1=t_0$ we see that $\Lambda_J$ exists for \emph{some} $J$.  Further, if for some collection $J_1,J_2,\ldots$ the $\Lambda_{J_i}$ exist, then $\Lambda_{J_i}|Y_1\times (J_i\cap J_j)=\Lambda_{J_j}|Y_1\times (J_i\cap J_j)$ by Theorem~\ref{UniquenessTheorem} part~\ref{UniquenessOne}, and so the $\Lambda_{J_i}$ can be patched to obtain $\Lambda_J$ where $J=\bigcup J_i$.

The last remark means that there is a \emph{largest} open interval $J_\mathrm{max}$ for which $\Lambda_{J_\mathrm{max}}$ exists.  We claim that $J_\mathrm{max}=I$.  If not then without loss there is a minimal $i\in I$ with $i>j$ for all $j\in J_\mathrm{max}$.  We may now take an open interval $J$, contained in $I$ and containing $i$, on which the conclusion of Lemma~\ref{UniquenessLemmaTwo} holds.  Take $t_1\in J\cap J_\mathrm{max}$, so that $\Lambda_{J_\mathrm{max}}|Y_1\times\{t_1\}$ extends to an integral $2$-morphism $L\colon\Phi|Y_1\times J\Rightarrow\Psi|Y_1\times J$ that, by Theorem~\ref{UniquenessTheorem} part~\ref{UniquenessOne}, coincides with $\Lambda|Y_1\times J_\mathrm{max}$ on $J_\mathrm{max}\cap J$.  Thus $\Lambda_{J_\mathrm{max}}$ and $L$ can be patched to obtain $\Lambda_{J\cup J_\mathrm{max}}$, contradicting the maximality of $J_\mathrm{max}$.  Thus $J_\mathrm{max}=I$ as claimed.

We have shown that for any $Y_1\subset Y$, open with compact closure, there is an integral $2$-morphism $\Lambda_{Y_1}\colon\Phi|Y_1\times I\Rightarrow\Psi|Y_1\times I$ extending $\lambda|Y_1$.  We may find a nested sequence of subsets $Y_1\subset Y_2\subset\cdots Y$ with compact closure and with $Y=\bigcup Y_i$, and we write $\Lambda_{Y_i}\colon\Phi|Y_i\times I\Rightarrow\Psi|Y_i\times I$ for the $2$-morphisms extending $\lambda|Y_i$ just obtained.  Then for any $i>j$, $\Lambda_{Y_i}|Y_j\times I=\Lambda_{Y_j}$ by Theorem~\ref{UniquenessTheorem} part~\ref{UniquenessOne}, and so the $\Lambda_{Y_i}$ can be patched to obtain the required integral $2$-morphism $\Lambda\colon\Phi\Rightarrow\Psi$ extending $\lambda$.
\end{proof}

We now move onto the proof of Theorem~\ref{ExistenceTheorem}.  In this result an assumption of properness has again been made.  In the proof of Theorem~\ref{UniquenessTheorem} the properness assumption was used to guarantee the lifting of integral curves over all times.  In the next proof, however, this assumption will be used to guarantee the existence of integral curves over some small time interval.  In what follows we will use the phrase \emph{at time $t$} to refer to what happens when one restricts a morphism or map $X\times I\to Y$, with $I$ an open interval, to the subset $X\times\{t\}$.  

\begin{proof}[Proof of Theorem~\ref{ExistenceTheorem}.]
Without loss we may assume that $\phi$ is the inclusion $\iota\colon\X_1\hookrightarrow\X$ of a full open substack with $\mathrm{cl}(\bar\X_1)\subset\bar\X$ compact.

Let $U\to\X$ be an atlas equipped with a vector field $X_U$ as in Corollary~\ref{VectorFieldCorollary}.  Now for each $x\in\mathrm{cl}(\bar\X_1)$ we may find an open subset $U_x\subset U$ that contains a representative of $x$ and has compact closure.  Since $\X$ is proper the map $U\times_\X U\to U\times U$ is proper, and so $U_x\times_\X U_x\subset U\times_\X U$ also has compact closure.  We may therefore find $\epsilon_x>0$ and maps
\begin{gather*}
\phi_x^0\colon U_x\times(-\epsilon_x,\epsilon_x)\to U\\
\phi_x^1\colon U_x\times_\X U_x\times(-\epsilon_x,\epsilon_x)\to U\times_\X U
\end{gather*}
which restrict at time $0$ to the inclusions and which integrate $X_U$, $X_{U\times_\X U}$ respectively.

Since $\mathrm{cl}(\bar\X_1)$ is compact we may find $x_1,\ldots,x_n\in\mathrm{cl}(\bar\X_1)$ such that each point of $\mathrm{cl}(\bar\X_1)$ is represented by a point of $\bigsqcup U_{x_i}$.  Setting $\tilde U=(\bigsqcup U_{x_i})_{\bar\X_1}$ and $\epsilon=\min\epsilon_{x_i}$, we obtain the following:
\begin{enumerate}
\item Maps $i_0\colon\tilde U\to U$, $i_1\colon\tilde U\times_\X\tilde U\to U\times_\X U$.
\item\label{SecondPart}A factorization
\[\xymatrix{
\tilde U\ar[r]^{i_0}\ar[d] & U\ar[d]_{}="1"\\
\X_1\ar@{^{(}->}[r]_\iota^{}="2"&\X\ar@{=}"1";"2"
}\]
in which $\tilde U\to\X_1$ is an atlas.  This diagram induces $i_1\colon\tilde U\times_\X\tilde U\to U\times_\X U$.
\item Maps $\phi_0\colon\tilde U\times(-\epsilon,\epsilon)\to U$, $\phi_1\colon\tilde U\times_\X\tilde U\times(-\epsilon,\epsilon)\to U\times_\X U$ integrating $X_U$, $X_{U\times_\X U}$ respectively and restricting to $i_0$, $i_1$ at time $0$.
\end{enumerate}
Write $\Upsilon$ for the Lie groupoid $U\times_\X U\rightrightarrows U$ representing $\X$ and write $\tilde\Upsilon$ for the groupoid $\tilde U\times_\X\tilde U\rightrightarrows\tilde U$ representing $\X_1$.  Then $i_0$ and $i_1$ form a groupoid morphism $i\colon\tilde\Upsilon\to\Upsilon$ and $X_U$ and $X_{U\times_\X U}$ form a groupoid morphism $X_\Upsilon\colon\Upsilon\to T\Upsilon$.  

Since $\phi_0$ and $\phi_1$ integrate $X_U$ and $X_{U\times_\X U}$ and restrict to $i_0$ and $i_1$ at time $0$, it is simple to verify that they define a groupoid map $\phi\colon\tilde\Upsilon\times(-\epsilon,\epsilon)\to\Upsilon$.  Dictionary Lemma~\ref{DictionaryLemma} part~\ref{DictionaryOne} then provides us with a morphism $\Phi\colon\X\times(-\epsilon,\epsilon)\to\X$ and a diagram
\[\xymatrix{
\tilde U\times(-\epsilon,\epsilon)\ar[r]^-{\phi_0}\ar[d]&U\ar[d]_{}="1"\\
\X_1\times(-\epsilon,\epsilon)\ar[r]_-\Phi^{}="2"&\X\ar@{=>}"1";"2"
}\]
that induces $\phi_1$.

We must prove that $\Phi$ is an integral of $X$ and that there is a $2$-morphism $\Phi|\X_1\times\{0\}\Rightarrow\iota$.  These are immediate consequences of Dictionary Lemma~\ref{DictionaryLemma} part~\ref{DictionaryTwo} .  First, since $X_\Upsilon\circ\phi=T\phi\circ \bigddt$, we obtain the required $2$-morphism $t_\Phi$; since $X_\Upsilon$ and $\bigddt$ are vector fields on the groupoids $\tilde\Upsilon\times(-\epsilon,\epsilon)$, $\Upsilon$, the condition on the resulting square \eqref{CategorifiedFlowDiagram} follows immediately.  Second, $\phi|\tilde\Upsilon\times\{0\}$ is just $i$, and so there is a $2$-morphism $\Phi|\X_1\times\{0\}\Rightarrow\iota$ as required.
\end{proof}

\begin{proof}[Proof of Theorem~\ref{RepresentabilityTheorem}.]
We shall prove that if $\Phi|\Y\times\{t_0\}$ is representable, then so is $\Phi$; the other direction is clear.  We claim that for any full open substack $\Y_1\subset\Y$ with $\mathrm{cl}(\bar\Y_1)\subset\bar\Y$ compact, and for any $t_1\in I$, there is an open interval $J\subset I$ containing $t_1$, with the property that for any $s\in J$, $\Phi|\Y_1\times J$ is representable if and only if $\Phi|\Y_1\times\{s\}$ is representable.

This claim allows us to prove the theorem.  For we can find $\epsilon>0$ such that $\Phi|\Y_1\times(t_0-\epsilon,t_0+\epsilon)$ is representable, so that if $\Phi|\Y_1\times I$ is not representable then we can without loss find a maximal $t_1\in I$ such that $\Phi|\Y_1\times(t_0-\epsilon,t_1)$ is representable.  Applying the claim again gives a contradiction, so that $\Phi|\Y_1\times I$ is representable.  Since $\Y$ is a nested union of such $\Y_1$, the theorem follows.

We now prove our claim.  Again, the `only if' part is trivial.  Let $U\to\X$ be an atlas and take $U_1\subset U$ open with compact closure such that $\Phi(\Y\times\{t_1\})$ is covered by $U_1$.  Choose $J$ so that $U_1\to\X$ extends to a submersion $U_1\times J\to\X\times I$ integrating $X\times\smallddt$.  By reducing $J$ if necessary we can assume that $U_1\times J\to\X\times I$ covers $\Phi\times\pi_2(\Y_1\times J)$.  Now $\Phi|\Y_1\times J$ is representable if and only if $\Y_1\times J\to\X\times I$ is representable, which is if and only if $(\Y\times J)\times_{\X\times I}(U_1\times J)$ is representable.  But $(\Y\times\{s\})\times_{\X}(U_1\times\{s\})$ is representable by assumption, so that $(\Y\times J)\times_{\X\times I}(U_1\times J)$ is representable by Lemma~\ref{PullbackProposition} below.  This completes the proof.
\end{proof}

\begin{lemma}[Pullbacks of integrals]\label{PullbackProposition}
Let $X$ be a vector field on a proper differentiable stack $\X$, and let $\Phi\colon\mathfrak{A}\times I\to\X$, $\Psi\colon\mathfrak{B}\times I\to\X$ be morphisms integrating $X$, with $\mathfrak{A}$ and $\mathfrak{B}$ differentiable, and further such that $\mathfrak{A}\times I\to\X\times I$ is a submersion.  Then the following diagram, whose $2$-morphism is furnished by part~\ref{UniquenessTwo} of Theorem~\ref{UniquenessTheorem}, is cartesian.
\begin{equation}\label{IntegralPullbackSquare}
\xymatrix{
(\mathfrak{A}\times\{0\})\times_\X(\mathfrak{B}\times\{0\})\times I\ar[r]\ar[d]&\mathfrak{A}\times I\ar[d]_{}="1"^{\Phi\times\pi_2}\\
\mathfrak{B}\times I\ar[r]_{\Psi\times\pi_2}^{}="2"&\X\times I\ar@{=>}"1";"2"
}
\end{equation}
\end{lemma}
\begin{proof}
We may assume that $\mathfrak{A}=A$ and $\mathfrak{B}=B$ are manifolds.  Then from the vector fields $\smallddt$ on $A\times I$, $\smallddt$ on $B\times I$, and $X\times\smallddt$ on $\X\times I$, we obtain a vector field $Y$ on $P=({A}\times I)\times_\X({B}\times I)\times I$ and a proper map $P\to (A\times I)\times (B\times I)$ that intertwines $Y$ and $\smallddt\times\smallddt$.  Now Lemma~\ref{ProperLemma} shows that every integral curve of $Y$ through $(A\times\{t_0\})\times_\X(B\times\{t_0\})$ can be defined over the entire time interval $I$, and that every point of $P$ lies on one of these flow lines.  This defines the required diffeomorphism $P\cong (A\times\{t_0\})\times_\X(B\times\{t_0\})\times I$.
\end{proof}

\begin{proof}[Proof of Theorem~\ref{FlowTheorem}.]
The first two parts are immediate from Theorem~\ref{UniquenessTheorem} and Proposition~\ref{RepresentabilityTheorem}.  We prove the third part.

Write $\X=\mathfrak{A}\cup\mathfrak{B}$ as a union of full open substacks with $\mathfrak{A}$ corresponding to the complement of $\mathrm{supp}(X)$ and with $\bar{\mathfrak{B}}$ containing $\mathrm{supp}(X)$ and such that $\mathrm{cl}(\bar{\mathfrak{B}})\subset\bar\X$ is compact.  Write $\iota_\mathfrak{A}$, $\iota_\mathfrak{B}$ for the inclusions.  Then by Theorem~\ref{ExistenceTheorem} we may find an open interval $J\supset\{0\}$ and a morphism $\phi_\mathfrak{B}\colon\mathfrak{B}\times J\to\X$ integrating $\X$ and admitting $\phi|\mathfrak{B}\times\{0\}\Rightarrow\iota_\mathfrak{B}$.  Since $X|\mathfrak{A}$ is equivalent to the zero-section we can find $\phi_\mathfrak{A}\colon\mathfrak{A}\times J\to\X$ integrating $X$ and admitting $\phi_\mathfrak{A}|\mathfrak{A}\times\{0\}\Rightarrow\iota_\mathfrak{A}$.

Theorem~\ref{UniquenessTheorem} allows us to patch $\phi_\mathfrak{A}$ and $\phi_\mathfrak{B}$ and obtain $\phi\colon\X\times J\to\X$ integrating $\X$ and admitting $\phi|\X\times\{0\}\Rightarrow\Id_\X$.  We will call such morphisms \emph{partial flows}.

Let $I\subset\R$ be the union of all those $J$ for which there exists a partial flow $\X\times J\to\X$.  Write $I=\bigcup_{i=1}^\infty J_i$ as a countable union of intervals for which there exist partial flows $\phi_i\colon\X\times J_i\to\X$.  Theorem~\ref{UniquenessTheorem} allows us to obtain a partial flow $\X\times I\to\X$.  $I$ is therefore the unique largest interval for which there is a partial flow $\phi_I\colon\X\times I\to\X$.  We claim that in fact $I=\R$.

If $I\neq\R$, then without loss $I$ is bounded above, so choose any positive $t_0\in I$.  Then both $\phi_I$ and
\[\X\times(I+t_0)\xrightarrow{(\phi_I\X\times\{t_0\})\times-t_0}\X\times I\xrightarrow{\phi_I}\X\]
integrate $X$ and are $2$-isomorphic when restricted to $\X\times\{t_0\}$, so by Theorem~\ref{UniquenessTheorem} can be glued to obtain a partial flow $\X\times I\cup(I+t_0)\to\X$ of $X$.  This contradicts the maximality of $I$.  Consequently $I=\R$ and the theorem is proved.
\end{proof}

\section{Global quotients}\label{QuotientsSection}
Let $G$ be a compact Lie group acting smoothly on a manifold $M$ and write $[M/G]$ for the quotient stack.  One expects that the geometry of $[M/G]$ is simply the $G$-equivariant geometry of $M$.  The results of this section are a clear instance of this principle, for we will show that $\Vect[M/G]$ can be described in terms of the $G$-invariant vector fields on $M$, and that the flow of a vector field on $[M/G]$ can be described in terms of the flow of the corresponding $G$-invariant vector field on $M$.

\begin{proposition}\label{EquivariantVectorFieldsProposition}
The groupoid $\Vect[M/G]$ of vector fields on $[M/G]$ is equivalent to the groupoid whose:
\begin{itemize}
\item objects are the $G$-invariant vector fields on $M$; 
\item arrows $X\to X'$ are the functions $\psi\colon M\to\mathfrak{g}$ such that $X'(m)=X(m)+\iota\psi(m)$ and $\psi(mg)=\mathrm{Ad}_{g^{-1}}\psi(m)$.
\end{itemize}
Here $\iota(v)$ denotes the tangent vector obtained by differentiating the $G$-action in the direction $v$.  The above equivalence restricts to an equivalence between the full subgroupoid on the compactly-supported invariant vector fields on $M$ and the full subgroupoid of compactly-supported vector fields on $[M/G]$.
\end{proposition}

\begin{proposition}\label{EquivariantFlowProposition}
Let $X$ be a compactly-supported vector field on $[M/G]$ corresponding to a $G$-invariant vector field $X_M$ on $M$ and let $\phi\colon M\times\R\to M$ be the flow of $X_M$.  The morphism of stacks
\[\Phi\colon[M/G]\times\R\to[M/G]\]
determined by $\phi$ is a flow of $X$.
\end{proposition}

\begin{proof}[Proof of Proposition~\ref{EquivariantVectorFieldsProposition}.]
Write $\Vect(M/G)$ for the groupoid described in the statement of the proposition.  We wish to find an equivalence $\Vect[M/G]\simeq\Vect(M/G)$.  Since the stack $[M/G]$ is represented by the action groupoid $M\rtimes G= M\times G\rightrightarrows M$, Theorem~\ref{VectorFieldTheorem} provides us with an equivalence $\Vect[M/G]\simeq\Vect(M\rtimes G)$.  Writing out $\Vect(M\rtimes G)$ explicitly (Definition~\ref{LieGpdVectorFieldDefinition}) and using $TG\cong G\times\mathfrak{g}$, we have an equivalence between $\Vect[M/G]$ and the groupoid whose:
\begin{itemize}
\item objects are pairs $(X,Y)$ consisting of a vector field $X$ on $M$ and a map $Y\colon M\times G\to\mathfrak{g}$ such that $X(mg)=X(m)g+\iota Y(m,g)$ and $Y(m,gh)=\mathrm{Ad}_{h^{-1}}Y(m,g)+Y(mg,h)$;
\item morphisms $(X,Y)\to(X',Y')$ are maps $\psi\colon M\to\mathfrak{g}$ for which $X'(m)=X(m)+\iota\psi(m)$ and $Y(m,g)+\psi(mg)=\mathrm{Ad}_{g^{-1}}\psi(m)+Y'(m,g)$.
\end{itemize}
It is clear from this description that $\Vect(M/G)$ is the full subgroupoid of $\Vect(M\rtimes G)$ on those objects $(X,Y)$ for which $Y=0$.  We will prove the claim $\Vect[M/G]\simeq\Vect(M/G)$ by showing that every object of $\Vect(M\rtimes G)$ is isomorphic to an object of $\Vect(M/G)$.

Fix a smooth invariant measure on $G$.  Let $(X,Y)$ be an object of $\Vect(M\rtimes G)$.  Define a vector field $\tilde X$ on $M$ by
\[\tilde X(m)=\int_{g\in G}X(mg)g^{-1}\]
and define $\psi\colon M\to\mathfrak{g}$ by
\[\psi(m)=\int_{g\in G}\mathrm{Ad}_{g}Y(m,g).\]
It is now routine to check that $(\tilde X,0)$ is an object of $\Vect(M/G)$ and that $\psi\colon(X,Y)\to(\tilde X,0)$ in $\Vect(M\rtimes G)$.

It remains to prove the second claim regarding the compactly-supported vector fields on $M$.  It is clear that a compactly-supported object of $\Vect(M/G)$ leads to a compactly supported vector field on $[M/G]$.  Conversely, take a compactly-supported vector field on $[M/G]$. This can, by the proof of Proposition~\ref{PatchEquivalenceProposition}, be represented by a vector field on $M\rtimes G$ that is not equal to the zero section only on a subgroupoid of $M\rtimes G$ whose image in $M/G$ has compact closure.  That is to say, the vector field on $M\rtimes G$ is given by compactly-supported vector fields on $M$ and $M\times G$.  The averaging process above clearly preserves this property, and the result follows.
\end{proof}

\begin{proof}[Proof of Proposition~\ref{EquivariantFlowProposition}.]
There is an obvious $2$-morphism $\Phi|[M/G]\times\{0\}\Rightarrow\Id_{[M/G]}$, and so it only remains to check that $\Phi$ integrates $X$.  But each arrow in the diagram
\[\xymatrix{
T[M/G]\times\R\ar[r]^{T\Phi}_{}="2" & T[M/G]\\
[M/G]\times\R\ar[u]^{\bigddt}\ar[r]_\Phi&[M/G]\ar[u]_X^{}="1"\ar@{.>}"1";"2"
}\]
is represented by a specific morphism of groupoids, constructed from $\phi$ or from $X_M$, and the corresponding diagram of groupoid morphisms commutes on the nose.  Thus, by Dictionary Lemma~\ref{DictionaryLemma} part~\ref{DictionaryTwo}, we may fill the square above with the required $2$-morphism $t_\Phi$.  The condition on \eqref{CategorifiedFlowDiagram} follows, again using the Dictionary Lemma.
\end{proof}

\section{\'Etale stacks}\label{EtaleSection}

A map $f\colon M\to N$ is \emph{\'etale} if it is a local diffeomorphism, or equivalently if the derivatives $T_m f$ are all linear isomorphisms.

\begin{definition}
A differentiable stack is \emph{\'etale} if it admits an \'etale atlas.
\end{definition}

An atlas $X\to\X$ is always representable, so it makes sense to ask whether it is also \'etale.  Note that $X\to\X$ is \'etale if and only if one or equivalently both of the projections $X\times_\X X\to X$ is \'etale.  Thus $\X$ is \'etale if and only if it is represented by a Lie groupoid whose source and target maps are \'etale.

Vector fields on manifolds have a particular functoriality under \'etale maps that they do not enjoy under general maps: they can be {pulled back}.  If $f\colon U\to V$ is \'etale and $X$ is a vector field on $V$ then the \emph{pullback} $f^\ast X$ denotes the vector field on $U$ given by $f^\ast X(u)=(T_uf)^{-1}X(f(u))$.  Note that $g^\ast f^\ast X=(fg)^\ast X$.

This functoriality of vector fields allowed the author in \cite{\Morse} to define vector fields and integral morphisms for \'etale stacks by considering the collection of all \'etale morphisms into the stack being studied.  To be precise:

\begin{definition}[{\cite[5.2]{\Morse}}]\label{EtaleVectorFieldDefinition}
A \emph{vector field} on an \'etale stack $\X$ is an assignment
\[ (U\to\X)\mapsto X_U\]
that sends each \'etale morphism from $U$ into $\X$ to a vector field on $U$.  This assignment is required to satisfy $f^\ast X_U=X_V$ whenever one has a triangle of \'etale morphisms
\[\xymatrix{
V\ar[rr]_{}="1"\ar[d]_f &{}&\X.\\
U\ar[rru]^{}="2"\ar@{=>}"1";"2" &{}&{}
}\]
\end{definition}

\begin{definition}[{\cite[5.7]{\Morse}}]\label{EtaleIntegralDefinition}
Let $X$ be a vector field on an \'etale stack $\X$.  Given a representable morphism $\Phi\colon\Y\times I\to\X$ and an \'etale morphism $U\to\X$, the pullback $(\Y\times I)\times_\X U$ is a manifold and its projection to $\Y\times I$ is \'etale.  Write $\Phi_U\colon (\Y\times I)\times_\X U\to U$ for the second projection.  We say that $\Phi$ \emph{is an integral morphism} if for each $U\to\X$ \'etale we have
\[T\Phi_U\left(\bigddt\right)=X_U\circ\Phi_U.\]
\end{definition}

These definitions are arguably more concrete and accessible than the definitions given in \S\ref{VectorFieldsSection} and \S\ref{IntegralsFlowsSection}.  They are certainly simpler in the sense that they do not require us to construct the tangent stack functor $T$.  In this section we are going to show how the two concepts above are equivalent to the ones established earlier.

\begin{lemma}
Let $U\to\X$ be \'etale.  Then the diagram
\begin{equation}\label{EtaleCartesianSquare}\xymatrix{
TU\ar[r]\ar[d] & U\ar[d]_{}="1"\\
T\X\ar[r]^{}="2"&\X \ar@{=>}"1";"2"
}\end{equation}
is cartesian.
\end{lemma}
\begin{proof}
Choose equivalences $\X\simeq\B(X_1\rightrightarrows X_0)$, $U\simeq\B(U_1\rightrightarrows U_0)$ where $X_0\to\X$, $U_0\to U$ are \'etale atlases and $U\to\X$ is obtained from a Lie groupoid morphism $(U_1\rightrightarrows U_0)\to(X_1\rightrightarrows X_0)$ that is \'etale in each component.  Then $U\times_\X T\X$ is represented by the Lie groupoid
\[TX_1\times_{X_0} X_1\times_{X_0}U_1\rightrightarrows TX_0\times_{X_0} X_1\times_{X_0}U_0\]
which, since all the maps forming the pullbacks are \'etale, is isomorphic to
\[T^\gpd (X_1\times_{X_0}X_1\times_{X_0}U_1\rightrightarrows X_0\times_{X_0} X_1\times_{X_0}U_0)\]
which is equivalent to $T^\gpd(U_1\rightrightarrows U_0)$ and, finally, $TU$.
\end{proof}

\begin{corollary}\label{EtaleBundleProposition}
If $\X$ is an \'etale stack then the projection $\pi_\X\colon T\X\to\X$ is a vector-bundle.
\end{corollary}
\begin{remark}
Corollary~\ref{EtaleBundleProposition} is in strong contrast to the general situation, in which the fibres of $\pi_\X$ can have the form $[V/W]$, where $V$ and $W$ are vector spaces and $W$ acts linearly on $V$.
\end{remark}

\begin{corollary}\label{EtaleVectorFieldCorollary}
Let $X$ be a vector field on an \'etale stack $\X$ and let $U\to\X$ be \'etale.  Then there is a \emph{unique} diagram
\[\xymatrix{
U\ar[r]^{X_U}\ar[d] & TU\ar[d]_{}="1"\\
\X\ar[r]_X^{}="2"&T\X\ar@{=>}"1";"2"
}\]
satisfying the conclusion of Corollary~\ref{VectorFieldCorollary}, and this diagram is cartesian.
\end{corollary}
\begin{proof}
Form the pullback of $TU\to T\X$ along $X$.  The fact that \eqref{EtaleCartesianSquare} is cartesian, together with the morphism $a_X\colon \pi_\X\circ X\Rightarrow\Id_\X$, identifies this pullback as $U$ and the resulting cartesian square has the form required by Corollary~\ref{VectorFieldCorollary}.  Any other diagram satisfying the conditions of Corollary~\ref{VectorFieldCorollary} is then related to this one by a map $U\to U$ which is necessarily the identity, and the diagrams therefore coincide.
\end{proof}

\begin{proposition}
Let $\X$ be an \'etale stack and let $\Vect^\mathrm{et}(\X)$ denote the set of vector fields on $\X$ as defined in Definition~\ref{EtaleVectorFieldDefinition}.  Regard $\Vect^\mathrm{et}(\X)$ as a groupoid with only identity arrows.  Then there is an equivalence
\begin{gather*}
\Vect(\X)\to\Vect^\mathrm{et}(\X)\\
X\mapsto ((U\to\X)\mapsto X_U)
\end{gather*}
where each $X_U$ is determined by Corollary~\ref{EtaleVectorFieldCorollary}.
\end{proposition}
\begin{proof}
Let $X$ be a vector field on $\X$.  Corollary~\ref{EtaleVectorFieldCorollary} provides us with vector fields $X_U$ on $U$ for each \'etale $U\to\X$; it further shows that the resulting assignment $(U\to\X)\mapsto X_U$ satisfies the conditions of Definition~\ref{EtaleVectorFieldDefinition}.  We therefore have a map from the objects of $\Vect(\X)$ to $\Vect^\mathrm{et}(\X)$.  But Corollary~\ref{EtaleVectorFieldCorollary} shows that if $X,Y$ are equivalent vector fields on $\X$ then their images in $\Vect^\mathrm{et}(\X)$ coincide.  Thus $\Vect(\X)\to\Vect^\mathrm{et}(\X)$ is a functor.

We now show that $\Vect(\X)\to\Vect^\mathrm{et}(\X)$ is fully faithful.
If $X\Rightarrow Y$ is an equivalence of vector fields on $\X$ then Corollary~\ref{EtaleVectorFieldCorollary} shows that the restrictions $X|U\Rightarrow Y|U$ are uniquely determined, and therefore such an equivalence, if it exists, is unique.  Moreover, if vector fields $X$ and $Y$ on $\X$ determine the same element of $\Vect^\mathrm{et}(\X)$ then Corollary~\ref{EtaleVectorFieldCorollary} determines a $2$-morphism $X|U\Rightarrow Y|U$ for each $U\to\X$ \'etale, and these satisfy the conditions required to ensure that they descend to an equivalence $X\to Y$.  This shows that $\Vect(\X)\to\Vect^\mathrm{et}(\X)$ is fully faithful.

We complete the proof by showing that $\Vect(\X)\to\Vect^\mathrm{et}(\X)$ is essentially surjective.  Any element $\{X_U\}$ of $\Vect^\mathrm{et}(\X)$ determines a morphism $X\colon\X\to T\X$ by choosing $U\to\X$ to be an \'etale atlas and considering the corresponding Lie groupoid.  That $X$ is a vector field determining the original $\{X_U\}$ is an immediate consequence of its construction.
\end{proof}

\begin{proposition}
Let $X$ be a vector field on an \'etale stack $\X$ and let $\Phi\colon\Y\times I\to\X$ be a representable morphism.  Then $\Phi$ integrates $X$ if and only if it satisfies the condition of Definition~\ref{EtaleIntegralDefinition}.
\end{proposition}
\begin{proof}
Suppose that $\Phi$ integrates $X$, let $U\to\X$ be \'etale, write $V\to\Y\times I$ for the induced \'etale atlas of $\Y\times I$, and let $\tilde\Phi\colon V\to U$ for the induced map; we are in the situation of diagram \eqref{SetupDiagramTwo}.  We must show that the diagram
\begin{equation}\label{CommutingSquare}\xymatrix{
V\ar[r]^{\tilde\Phi}\ar[d]_{\bigddt}& U\ar[d]^{X_U}\\
TV\ar[r]_{T\tilde\Phi}& TU
}\end{equation}
commutes.

By its construction the composite $V\to U\to TU$ fits into the $2$-commutative rectangle
\[\xymatrix{
V\ar[r]\ar[d]&U\ar[r]^{X_U}\ar[d]_{}="1"& TU\ar[d]_{}="2"\ar[r] & U\ar[d]_{}="3"\\
\Y\times I\ar[r]^{}="4"&\X\ar[r]^{}="5"&T\X\ar[r]^{}="6" &\X.\ar@{=>}"1";"4"\ar@{=>}"2";"5"\ar@{=>}"3";"6"
}\]
whose middle square is obtained using Corollary~\ref{EtaleVectorFieldCorollary}.  Since $TU\simeq U\times_\X T\X$, this rectangle determines the composition $V\to U\to TU$.  Similarly, $V\to TV\to TU$ is determined by a rectangle
\[\xymatrix{
V\ar[r]\ar[d]&TV\ar[r]\ar[d]_{}="1"& TU\ar[d]_{}="2"\ar[r] & U\ar[d]_{}="3"\\
\Y\times I\ar[r]^{}="4"&T(\Y\times I)\ar[r]^{}="5"&T\X\ar[r]^{}="6" &\X.\ar@{=>}"1";"4"\ar@{=>}"2";"5"\ar@{=>}"3";"6"
}\]
where now the first square is determined by Corollary~\ref{EtaleVectorFieldCorollary}.

Now we could paste the first rectangle with $t_\Phi$, and compose the $2$-morphisms.  Using the conditions on $t_\Phi$ and Corollary~\ref{EtaleVectorFieldCorollary}, we find that the composed $2$-morphism is identical with the composition of $2$-morphisms in the second of the rectangles.  Thus the two compositions $V\to TU$ are related by a $2$-morphism, so that in fact they coincide.  This shows that the square \eqref{CommutingSquare} above does indeed commute.

Conversely, if for each $U\to\X$ \'etale the diagram \eqref{CommutingSquare} commutes, then in particular it commutes when $U\to\X$ is taken to be an atlas $X_0\to\X$ or either of the induced maps $X_1\to\X$, where $X_1=X_0\times_\X X_0$.  Thus we have a commuting square of Lie groupoids that represents the required commuting square \eqref{FlowDiagram}. By its construction this square satisfies the condition on \eqref{CategorifiedFlowDiagram}.  This completes the proof.
\end{proof}

\appendix

\section{The Dictionary Lemma}
The \emph{dictionary lemma} below explains how to relate morphisms and $2$-morphisms of Lie groupoids
\[\xymatrix{
\Gamma\ar@/^2ex/[rr]_{}="1"\ar@/_2ex/[rr]^{}="2" &{}& \Delta\ar@{=>}"1";"2"
}\]
with morphisms and $2$-morphisms of stacks
\[\xymatrix{
\B\Gamma\ar@/^2ex/[rr]_{}="1"\ar@/_2ex/[rr]^{}="2" &{}& \B\Delta.\ar@{=>}"1";"2"
}\]

\begin{lemma}[The dictionary lemma, {\cite[2.6]{\BX}}]\label{DictionaryLemma}\hfill
\begin{enumerate}
\item\label{DictionaryOne}
A groupoid morphism $f\colon\Gamma\to\Delta$ determines a diagram
\[\xymatrix{
\Gamma_0\ar[r]^{f_0}\ar[d] & \Delta_0\ar[d]_{}="1"\\
\B\Gamma\ar[r]_f^{}="2" & \B\Delta\ar@{=>}"1";"2"_\eta
}\]
for which the induced map $\Gamma_0\times_{\B\Gamma}\Gamma_0\to\Delta_0\times_{\B\Delta}\Delta_0$ is just $f_1\colon\Gamma_1\to\Delta_1$.
\item\label{DictionaryTwo}
If a second diagram
\[\xymatrix{
\Gamma_0\ar[r]^{f_0}\ar[d] & \Delta_0\ar[d]_{}="1"\\
\B\Gamma\ar[r]_{f'}^{}="2" & \B\Delta\ar@{=>}"1";"2"_{\eta'}
}\]
has the same property as the diagram in part \ref{DictionaryOne}, then there is a unique $\epsilon\colon f\Rightarrow f'$ such that $\epsilon|\Gamma_0\circ\eta=\eta'$.
\item\label{DictionaryThree}
Let $f\colon\Gamma\to\Delta$ be Lie groupoid morphisms and let
\[\xymatrix{
\Gamma_0\ar[r]^{f_0}\ar[d] & \Delta_0\ar[d]_{}="1"\\
\B\Gamma\ar[r]_f^{}="2" & \B\Delta\ar@{=>}"1";"2"_\eta
}\qquad
\xymatrix{
\Gamma_0\ar[r]^{g_0}\ar[d] & \Delta_0\ar[d]_{}="1"\\
\B\Gamma\ar[r]_g^{}="2" & \B\Delta\ar@{=>}"1";"2"_\mu
}\]
be diagrams satisfying the property of part \ref{DictionaryOne}.  Then any 2-morphism $\phi\colon f\Rightarrow g$ can be composed with these diagrams to obtain
\[\xymatrix{
\Gamma_0\ar[r]^{f_0}\ar[d]_{f_0} & \Delta_0\ar[d]_{}="1"\\
\Delta_0\ar[r]^{}="2" & \B\Delta,\ar@{=>}"1";"2"
}\]
i.e.~a map $\tilde\phi\colon\Gamma_0\to\Delta_1$.  This $\tilde\phi$ is in fact a groupoid $2$-morphism $\tilde\phi\colon f\to g$.  This process determines a correspondence between 2-morphisms $f\Rightarrow g$ and $2$-morphisms $f\Rightarrow g$.
\end{enumerate}
\end{lemma}

\begin{note}
Not all stack morphisms $\B\Gamma\to\B\Delta$ arise in this way from groupoid morphisms $\Gamma\to\Delta$.  One might first have to replace $\Gamma$ be some refinement $\tilde\Gamma$.  One popular solution to this problem is to enlarge the $2$-category of Lie groupoids to the weak $2$-category $\mathsf{Bi}$ whose morphisms are `bibundles' \cite{\Lerman}.
\end{note}

\section{Proper stacks}\label{ProperAppendix}

A smooth map $f\colon M\to N$ is called \emph{proper} if for any compact $K\subset N$ the preimage $f^{-1}K$ is also compact.  This property of maps is local on the base and stable under pullbacks.  If $\X$ is a differentiable stack then the diagonal morphism $\Delta\colon\X\to\X\times\X$ is representable and we are able to make the following definition.

\begin{definition}
A differentiable stack $\X$ is \emph{proper} if the diagonal morphism $\Delta\colon\X\to\X\times\X$ is proper.
\end{definition}

If $X\to\X$ is an atlas then $X\times X\to\X\times\X$ is also an atlas and $(X\times X)\times_{\X\times\X}\X$ is equivalent to $X\times_\X X$.  It follows that $\X$ is proper if and only if it is represented by a Lie groupoid $\Gamma$ whose diagonal $s\times t\colon\Gamma_1\to\Gamma_0\times\Gamma_0$ is proper.

The proper \'etale differentiable stacks, also called differentiable Deligne-Mumford stacks, are precisely the orbifolds, and already present a significantly richer collection of objects than just manifolds.  Global quotients and gerbes, however, are proper but usually not \'etale.

\begin{theorem*}[{Zung, \cite[Theorem 2.3]{\Zung}}]
A proper Lie groupoid $\Gamma$ with fixed point $m\in\Gamma_0$ is locally isomorphic to the action groupoid $T_m\Gamma_0/\mathrm{Aut}_m$.
\end{theorem*}

\begin{corollary}\label{LocalQuotientCorollary}
A proper differentiable stack locally has the form of a global quotient $[M/G]$ with $G$ compact.
\end{corollary}
\begin{proof}
Let $\X$ be a proper differentiable stack and fix a point in $\X$.  Let $X\to\X$ be an atlas and choose a point $x\in X$ that represents the chosen point of $\X$.  Consider the proper groupoid $X\times_\X X\rightrightarrows X$.  By \cite[2.2]{\Zung} we can find an embedded submanifold $U\hookrightarrow X$ that contains $x$ and is such that $x$ is a fixed point of $U\times_\X U\rightrightarrows U$.  Moreover, by reducing $U$ if necessary we may assume that $U\hookrightarrow X$ is everywhere transverse to the orbits of $X\times_\X X\rightrightarrows X$.  It follows that $U\to\X$ is a submersion.  Consider the open substack $\mathfrak{U}$ of $\X$ whose atlas is $U$; this contains the chosen point of $\X$, and so we can prove the corollary by showing that $\mathfrak{U}$ is a global quotient. $\mathfrak{U}$ is represented by the proper Lie groupoid $U\times_\X U\rightrightarrows U$, and by the theorem above we may reduce $U$ one last time and assume that it in fact has the form $T_xU\times\mathrm{Aut}_x\rightrightarrows T_x U$, so that $\mathfrak{U}\simeq [T_x U/\mathrm{Aut}_x]$ as required.
\end{proof}

Recall the notion of \emph{underlying space} or \emph{orbit space} $\bar\X$ of a differentiable stack $\X$.  See \cite[\S2]{\Morse}.  The underlying space is a topological space derived from $\X$ and which is naturally homeomorphic to the orbit space of any Lie groupoid representing $\X$.  Open subsets of $\bar\X$ correspond to the full open substacks of $\X$.

\begin{definition}\label{POUDefinition}
A differentiable stack $\X$ \emph{admits smooth partitions of unity} if for each open cover $\{U_\alpha\}$ of $\bar\X$ there is a countable family $\{\phi_i\}$ of morphisms $\X\to\R$ such that the maps $\bar\phi_i\colon\bar\X\to\R$ form a partition of unity subordinate to $\{U_\alpha\}$.
\end{definition}

\begin{proposition}\label{POUProposition}
A proper differentiable stack admits partitions of unity.
\end{proposition}
\begin{proof}
The proof of this result is an immediate generalization of the proofs of the analogues for proper \'etale stacks given in \cite[\S3]{\Morse}.  See also \cite[\S3]{\JeffJohannes}.
\end{proof}

\begin{question}
Is a differentiable stack that admits partitions of unity and that is locally isomorphic to a quotient $[M/G]$ with $G$ compact necessarily proper?
\end{question}

\bibliographystyle{alpha}
\bibliography{Bibliography}
\end{document}